\documentclass[final,leqno,onefignum,onetabnum]{siamltex1213}

\usepackage{tikz}
\usepackage{relsize}
\usepackage{graphicx}
\usepackage{subfig}
\usepackage{enumitem}
\usepackage{amsmath,amssymb}
\usepackage{caption}

\newcommand{\E}{\mathbb{E}}

\newcommand{\Samp}{\Theta}
\newcommand{\samp}{\theta}
\newcommand{\loss}{\bar{\omega}}

\newcommand{\bxi}{\boldsymbol\xi}
\newcommand{\bp}{\mathbf{p}}

\newcommand{\bs}{\boldsymbol{\sigma}}

\title{OPTIMIZATION OF GAUSSIAN RANDOM FIELDS} 

\author{Eric Dow\thanks{Department of Aeronautics and Astronautics, Massachusetts Institute of Technology, Cambridge, MA, 02139. 
(\email{ericdow@mit.edu}).} \and Qiqi Wang\thanks{Department of Aeronautics and Astronautics, Massachusetts Institute of Technology, Cambridge, MA, 02139. 
(\email{qiqi@mit.edu}).}}

\begin{document}
\maketitle
\slugger{mms}{xxxx}{xx}{x}{x--x}%slugger should be set to mms, siap, sicomp, sicon, sidma, sima, simax, sinum, siopt, sisc, or sirev

\begin{abstract}
Many engineering systems are subject to spatially distributed uncertainty, i.e. uncertainty that can be modeled as a random field. Altering the mean or covariance of this uncertainty will in general change the statistical distribution of the system outputs. We present an approach for computing the sensitivity of the statistics of system outputs with respect to the parameters describing the mean and covariance of the distributed uncertainty. This sensitivity information is then incorporated into a gradient-based optimizer to optimize the structure of the distributed uncertainty to achieve desired output statistics. This framework is applied to perform variance optimization for a model problem and to optimize the manufacturing tolerances of a gas turbine compressor blade.

%Novelty:
%To the author's knowledge, this is the first paper that considers optimizing spatially distributed uncertainty. This work extends methods from probabilistic sensitivity analysis to handle random fields, and presents novel methods for sensitivity analysis of Gaussian random fields. Novel applications to tolerance optimization are also presented.
\end{abstract}



\pagestyle{myheadings}
\thispagestyle{plain}
\markboth{E. DOW AND Q. WANG}{OPTIMIZATION OF GAUSSIAN RANDOM FIELDS}

\section{Introduction and motivation}

An engineering system maps a set of inputs to a set of outputs, which quantify the performance of the system. In a deterministic design setting, the inputs are assumed to take a single (nominal) value, and the resulting outputs are deterministic functions of the nominal input values. In many engineering systems, the inputs are subject to some uncertainty due to natural variations in the system's environment or due to a lack of knowledge. In this case, the inputs can be modeled as random variables, and the system outputs are also, in general, random variables. The system performance is commonly quantified in terms of the statistics of the outputs, e.g. their mean or variance. The statistical distribution of the system outputs can be changed by either changing the distribution of the input uncertainty, or by changing the design of the system, i.e. how the inputs are mapped to the outputs. Design under uncertainty, also referred to as robust design, is often applied to optimize systems with random outputs. Broadly speaking, robust design methodologies construct designs whose performance remains relatively unchanged when the inputs are perturbed from their nominal value as a result of uncertainty\cite{beyer}. Examples include topology optimization of structures subject to random field uncertainties, design of gas turbine compressor blades subject to manufacturing variations, and optimization of airfoils subject to geometric uncertainty\cite{chen_2010, garzon_thesis, schillings_2011}. In these works, the system design is optimized to minimize the impact of variability on the output statistics.

In most applications of robust optimization, the statistical distribution of the input variability is assumed to be constant. In some applications, however, the distribution of the input uncertainty can be controlled. A concrete example is a gas turbine compressor blade subject to geometric variability introduced by the manufacturing process. In this context, the system inputs include the geometry of the compressor blade, which is assumed to be random as a result of random perturbations introduced by the manufacturing process. As will be described in the next section, the randomness in the blade geometry is an example of spatially distributed uncertainty, and can therefore be modeled as a random field. The outputs are chosen to describe the aerothermal performance of the compressor blade, e.g. the total pressure loss coefficient and flow turning. The mean performance of manufactured compressor blades has been shown to degrade as the level of variability (quantified by its standard deviation) increases\cite{garzon_thesis}. The level of variability can be reduced by specifying stricter manufacturing tolerances. However, specifying stricter manufacturing tolerances incurs higher manufacturing costs. Therefore, the cost associated with reducing variability competes with the benefits of improving performance, implying that there may be some optimal level of variability that balances these competing costs. 

This paper presents a method for optimizing the statistical distribution of random fields that describe the variability in a system's inputs. An efficient approach for computing the sensitivity of system outputs with respect to the parameters defining the distribution of the random field is presented. This sensitivity information is then used by a gradient-based optimizer to optimize these parameters. We apply this framework to perform variance optimization for a model problem as well as to a compressor blade tolerance optimization problem.

\section{Gaussian random fields}

\begin{figure}[htbp]
\centering
\subfloat[Realizations of a smooth random field]{\includegraphics[width=0.49\textwidth]{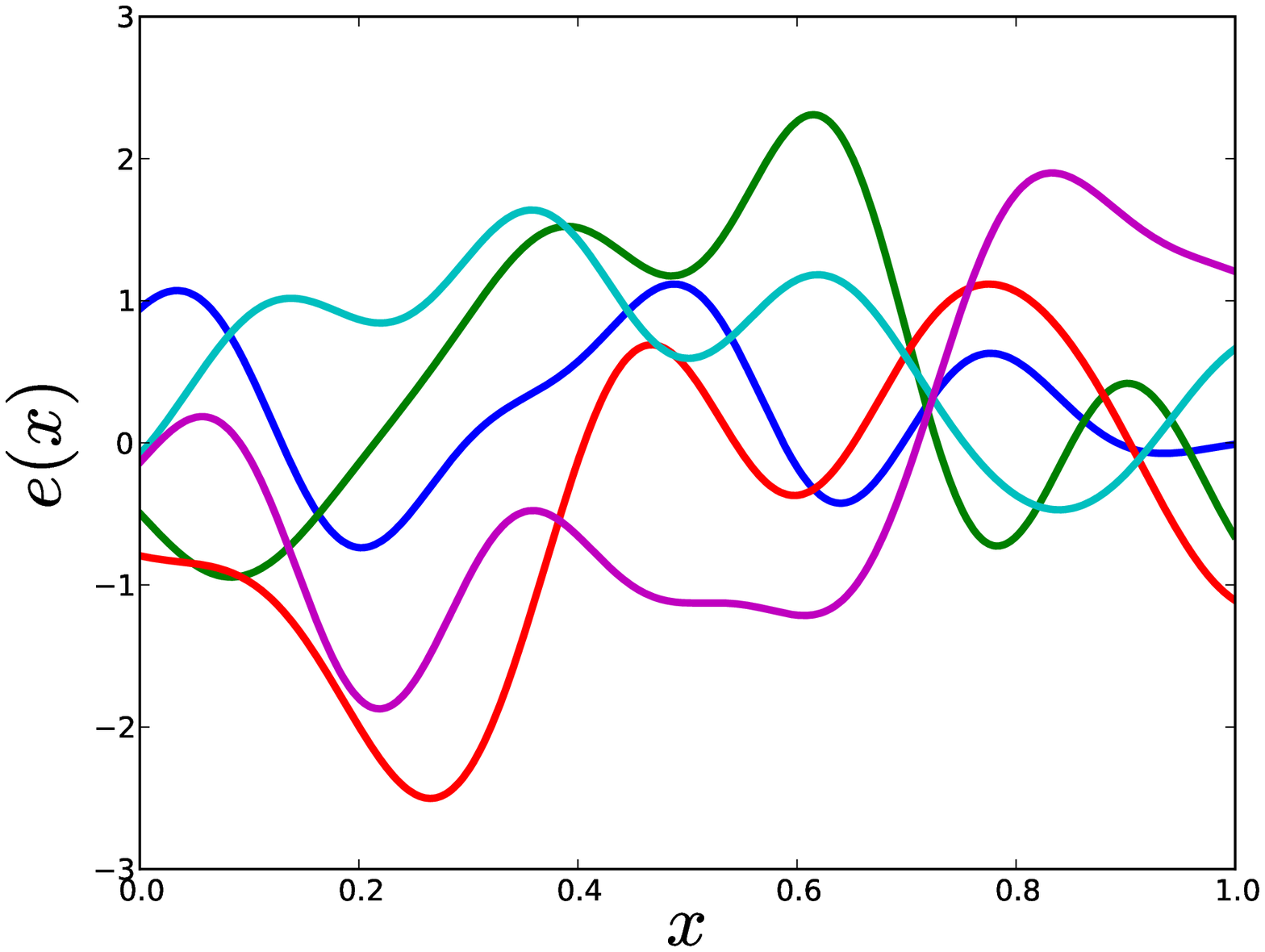}}
\hspace{0.1cm}
\subfloat[Realizations of a non-smooth random field]{\includegraphics[width=0.49\textwidth]{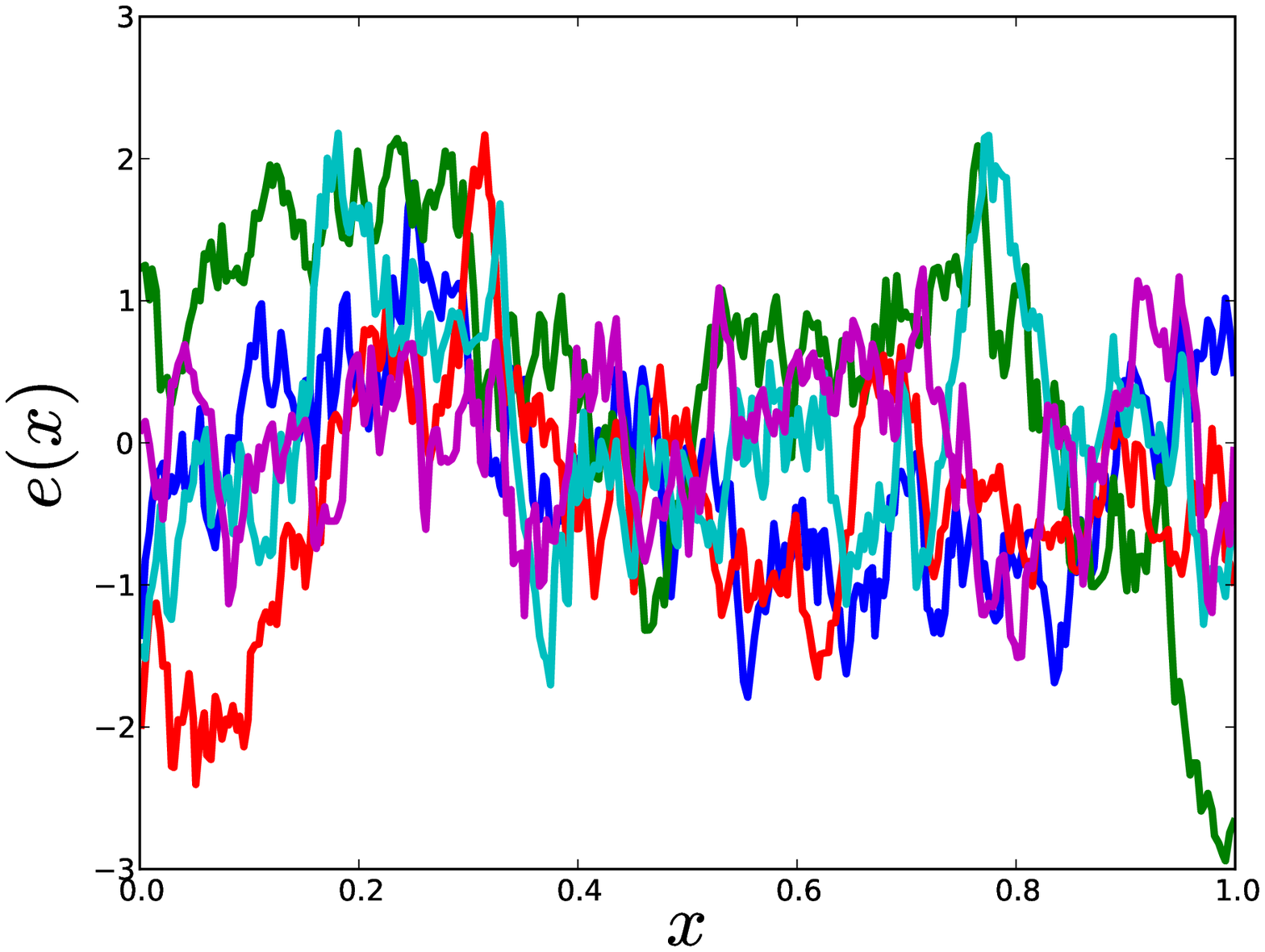}}
\newline
\subfloat[Realizations of a random field with large correlation length]{\includegraphics[width=0.49\textwidth]{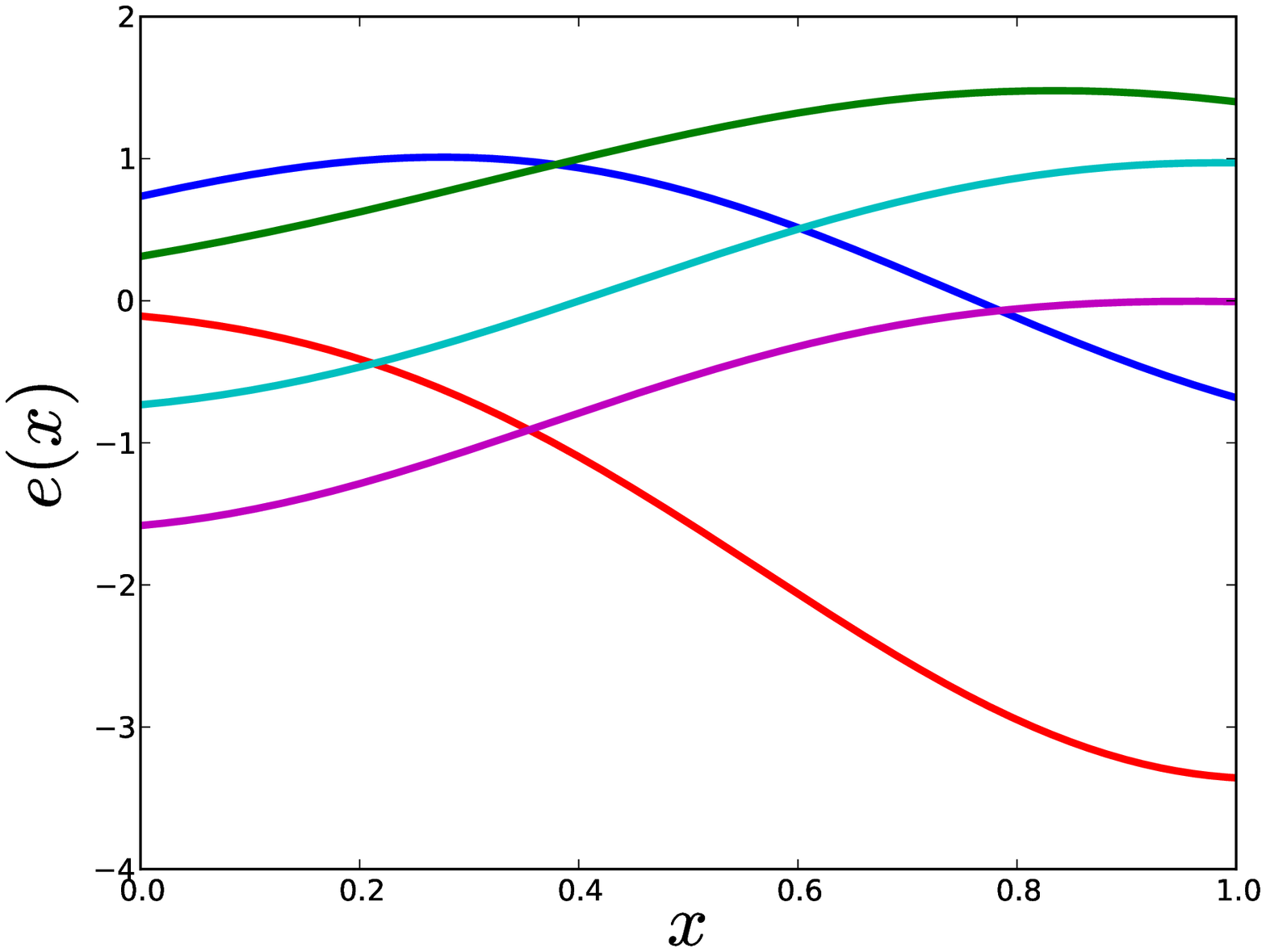}}
\hspace{0.1cm}
\subfloat[Realizations of a random field with a short correlation length]{\includegraphics[width=0.49\textwidth]{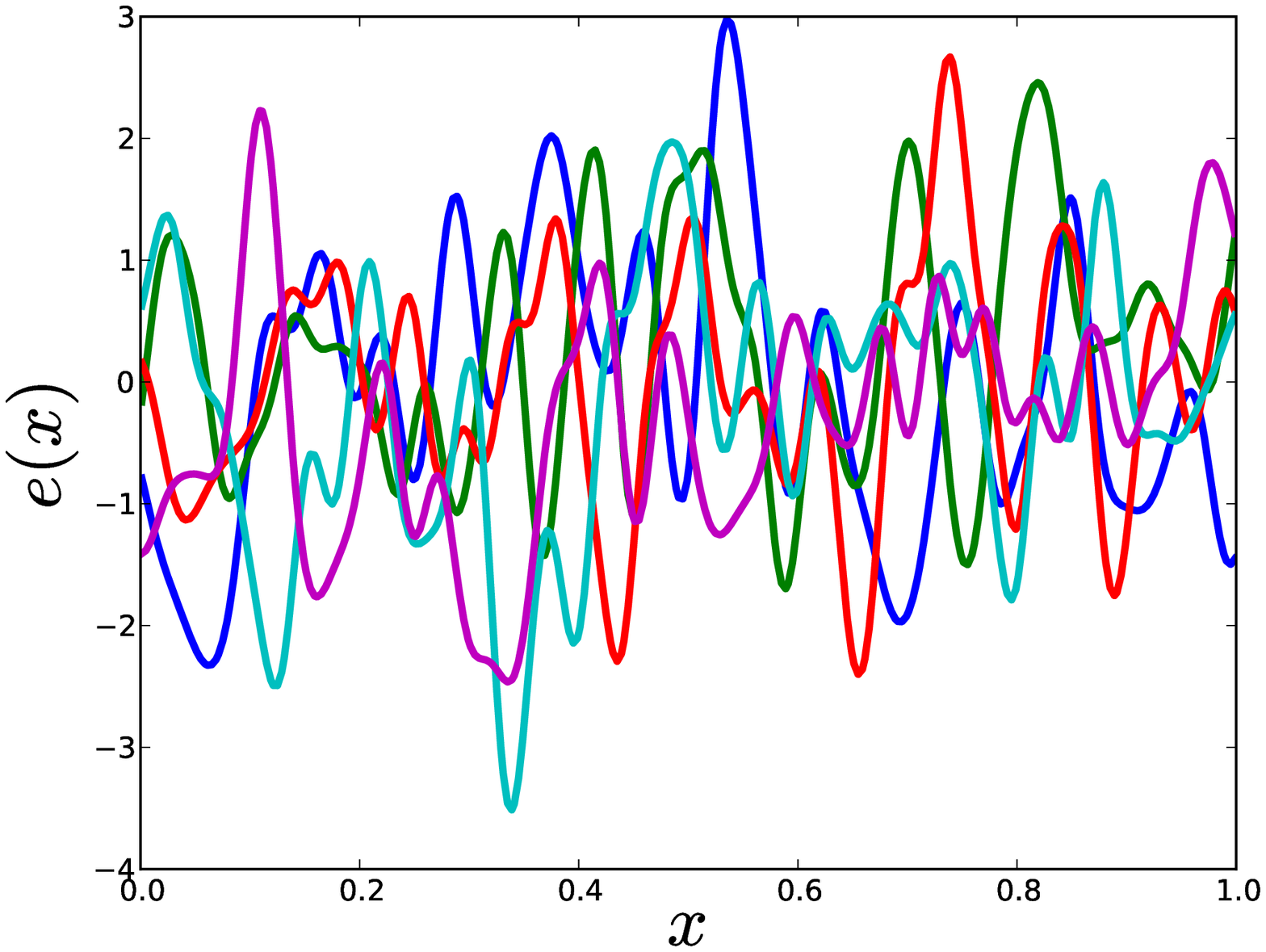}}
\caption{Illustration of the effect of covariance function on the smoothness and correlation length of a random field.}
\label{fig:rf_reals}
\end{figure}

Random fields provide a convenient method for modeling spatially distributed uncertainty. Random fields have previously been used to model spatially distributed uncertainty in a wide variety of systems, including natural variations in ground permeability\cite{christakos}, random deviations in material properties for structural optimization problems\cite{chen_2010}, and geometric variability in airfoils\cite{borzi_2010, schillings_2011}. Given a probability space $(\Samp, \mathcal{F}, \mathbb{P})$ and a metric space $X$, a random field is a measurable mapping $e : \Samp \rightarrow \mathbb{R}^X$\cite{adler}. In this work, we consider spatially distributed uncertainty in the form of a Gaussian random field $e(x,\samp)$. The defining characteristic of Gaussian random fields is that for any $x_1,...,x_n$, the vector $(e(x_1,\samp),...,e(x_n,\samp))$ is distributed as a multivariate Gaussian. Gaussian random fields are uniquely defined by their mean $\bar{e}(x)$ and covariance function $C(x_1,x_2)$:
\begin{equation}
\bar{e}(x) = \mathbb{E}[e(x,\samp)],
\end{equation}
\begin{equation}
C(x_1,x_2) = \mathbb{E}[(e(x_1,\samp) - \bar{e}(x_1))(e(x_2,\samp) - \bar{e}(x_2))],
\end{equation}
where the expectation is taken over $\samp$. The covariance function describes the smoothness and correlation length of the random field. Figure \ref{fig:rf_reals} shows realizations of random fields with different covariance functions. The realizations in the top left, produced with the squared exponential covariance function, are infinitely differentiable, and thus appear very smooth. Conversely, the realizations on the top right, produced with the exponential covariance function, are nowhere differentiable, and thus appear very jagged. The effects of changing the correlation length for the squared exponential kernel are shown in the bottom figures.

\subsection{The Karhunen-Lo\`{e}ve Expansion}

The Karhunen-Lo\`{e}ve (K-L) expansion, also referred to as the proper orthogonal decomposition (POD), can be used to represent a random field as a spectral decomposition of its covariance function\cite{lemaitre}. The random field $e$ is assumed to be continuous in the mean square sense:
\begin{equation}
\lim_{x_1 \rightarrow x_2}\E[(e(x_1,\samp) - e(x_2,\samp))^2] = 0\ \ \ \forall x_2 \in X.
\end{equation}
Then, the covariance function $C$ is continuous and 
\begin{equation}
\int_X \int_X C(x_1,x_2)\ dx_1dx_2 < \infty.
\end{equation}
We can therefore define the covariance kernel $K$ as
\begin{equation}
\langle Kv,w\rangle = \int_X \int_X C(x_1,x_2)v(x)w(x)\ dx_1\ dx_2,
\end{equation}
which is a symmetric semi-positive definite operator equipped with inner product $\langle\cdot,\cdot\rangle$, and $v$, $w \in L_2(X)$. By Mercer's theorem, it follows that $C$ has the spectral decomposition
\begin{equation}
C(x_1,x_2) = \sum_{i \geq 1} \lambda_i \phi_i(x_1) \phi_i(x_2),
\end{equation}
where each pair of eigenvalues $\lambda_i$ and eigenfuctions $\phi_i(x)$ are computed from the following Fredholm equation:
\begin{equation}
\int_S C(x_1,x_2)\phi_i(x_2)dx_2 = \lambda_i \phi_i(x_1).
\label{eq:fred}
\end{equation}
Moreover, the eigenfunctions can be chosen orthonormal such that $\langle \phi_i,\phi_j\rangle = \delta_{ij}$, and the eigenvalues are real, non-negative, and satisfy
\begin{equation}
\sum_{i \geq 1} \lambda_i^2 < \infty.
\end{equation}
By the Karhunen-Lo\`{e}ve theorem, the decomposition of the random field is given by:
\begin{equation}
e(x,\samp) = \bar{e}(x) + \sum_{i \geq 1} \sqrt{\lambda_i} \phi_i(x) \xi_i(\samp),
\label{eq:kl_full}
\end{equation}
where the eigenvalues are arranged in descending order such that $\lambda_1 \geq \lambda_2 \geq ... \rightarrow 0$. The distribution of the random variables $\xi_i(\samp)$ can be determined by taking the inner product of the random field with each of the eigenfunctions:
\begin{equation}
\xi_i(\samp) = \frac{1}{\sqrt{\lambda_i}}\langle e(x,\samp) - \bar{e}(x),\phi_i(x)\rangle.
\end{equation}
The random variables $\xi_i(\samp)$ are mutually uncorrelated with zero mean and unit variance.
%\begin{equation}
%\mathbb{E}[\xi_i] = \mathbb{E}\left[\frac{1}{\sqrt{\lambda_i}}\langle e - \bar{e},\phi_i(x)\rangle \right] 
%= \frac{1}{\sqrt{\lambda_i}}\langle\mathbb{E}[e - \bar{e}],\phi_i(x)\rangle = 0.
%\end{equation}
%\begin{eqnarray}
%\mathbb{E}[\xi_i\xi_j] & = & \mathbb{E}\left[\frac{1}{\sqrt{\lambda_i \lambda_j}} \int_X \int_X (e(x_1,\samp) - \bar{e}(x_1)) (e(x_2,\samp) - \bar{e}(x_2))\phi_i(x_1) \phi_j(x_2)\ dx_1 dx_2 \right] \notag \\
%& = & \frac{1}{\sqrt{\lambda_i\lambda_j}} \int_X \int_X \mathbb{E} \left[ (e(x_1,\samp) - \bar{e}(x_1)) (e(x_2,\samp) - \bar{e}(x_2))\right]\phi_i(x_1) \phi_j(x_2)\ dx_1 dx_2 \notag \\
%& = & \frac{1}{\sqrt{\lambda_i\lambda_j}} \int_X \left[\int_X C(x_1,x_2)\phi_i(x_1) \ dx_1 \right]\phi_j(x_2)\ dx_2 \notag \\
%& = & \frac{\lambda_i}{\sqrt{\lambda_i\lambda_j}} \langle \phi_i,\phi_j \rangle \notag \\
%& = & \delta_{ij} \notag
%\end{eqnarray}
For a Gaussian random field, the $\xi_i(\samp)$ are independent, identically distributed (i.i.d.) standard normal random variables.

To construct the K-L expansion numerically, the Nystr\"{o}m method is used\cite{nystrom}. The domain $X$ is discretized, and quadrature is used to approximate Equation (\ref{eq:fred}).
This results in a discrete eigenproblem of the form
\begin{equation}
\mathbf{C}\phi_i = \lambda_i \phi_i,
\label{eq:eig}
\end{equation}
where $\mathbf{C}$ is the discretized covariance matrix. Solving this eigenproblem gives the eigenvalues and eigenvectors evaluated on the discretized domain. The K-L expansion (\ref{eq:kl_full}) is truncated at a finite number of terms, resulting in an approximate spectral expansion of the random field:
\begin{equation}
\hat{e}(x,\samp) = \bar{e}(x) + \sum_{i=1}^{N_{KL}}\sqrt{\lambda_i}\phi_i(x)\xi_i(\samp).
\label{eq:kl}
\end{equation}
The truncated expansion minimizes the mean square error, and the decay of the eigenvalues determines the rate of convergence. The level of truncation $N_{KL}$ is often set equal to the smallest $k$ such that the partial scatter $S_k$ exceeds some threshold, where the partial scatter is defined as
\begin{equation}
S_k = \frac{\sum_{i=1}^k \lambda_i}{\sum_{i=1}^{N_s} \lambda_i}.
\end{equation}

\section{Optimizing the mean and covariance}

Consider a system whose performance is subject to spatially distributed uncertainty in the form of a Gaussian random field $e(x,\samp)$. Each output of the system is a functional of this random field, i.e. $F(\samp) = F(e(x,\samp))$, and is itself a random variable. $F$ can either be a direct functional of the random field, or a functional of the solution of a system of equations subject to random field uncertainty, e.g. the Navier-Stokes equations on a domain with a boundary that is described by $e(x,\samp)$. We are interested in the statistics $s_F$ of this functional, e.g. its mean or variance. In the case of multiple system output statistics, we generalize to the vector of output statistics $\mathbf{s}_\mathbf{F}$.

\begin{figure}[htbp]
\begin{center}
\includegraphics[scale=0.9]{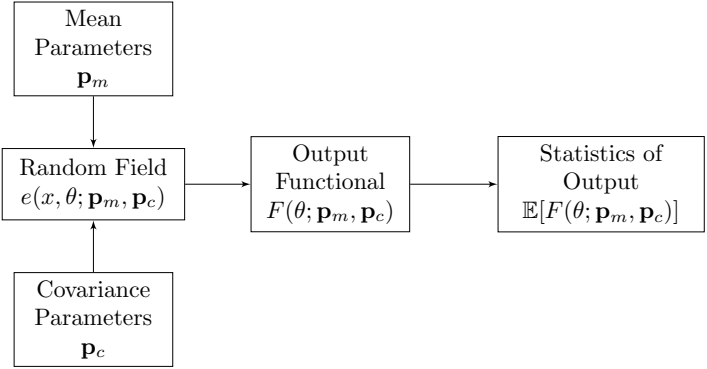}
\end{center}
\caption{Propagation of distributed uncertainty to the statistics of an output quantity of interest (in this case, the mean of the functional $F$).}
\label{fig:block}
\end{figure}

We aim to optimize the system's statistical response $s_F$ by controlling the mean and covariance of the random field $e(x,\samp)$. The design variables are then the mean of the random field $\bar{e}(x)$, parameterized by the vector $\bp_m$, and covariance of the random field $C(x_1,x_2)$, parameterized by the vector $\bp_c$. The design vector $\bp = \{\bp_m,\bp_c\}$ fully defines the Gaussian random field. We assume that $\bar{e}(x;\bp_m)$ and $C(x_1,x_2;\bp_c)$ depend smoothly on $\bp_m$ and $\bp_c$, respectively. Changing the mean and covariance of the random field will in general change the system output statistics, so that $s_F = s_F(\bp)$. Figure \ref{fig:block} illustrates the propagation of the random field to the output statistic $s_F = \mathbb{E}[F]$.

To optimize the statistical response of the system, we formulate the following optimization problem:
\begin{equation}
\begin{aligned}
& \bp^* =
& & \underset{\bp \in \mathcal{P}}{\arg \min}
& & f(\bp,\mathbf{s}_\mathbf{F}(\bp)) \\
& & & \ \ \ \ \text{s.t.}
& &  g(\bp,\mathbf{s}_\mathbf{F}(\bp)) \leq 0 \\
& & & & & h(\bp,\mathbf{s}_\mathbf{F}(\bp)) = 0,
\end{aligned}
\label{eq:gen_opt}
\end{equation}
where the objective and constraint functions $f$, $g$, and $h$ may depend on both the design parameters $\bp$ and the system output statistics $\mathbf{s}_\mathbf{F}(\bp)$, and $\mathcal{P}$ is the design space for the mean and covariance parameters. Note that, in general, the objective and constraint functions are nonlinear with respect to $\bp$.

\section{Sample average approximation}

To solve (\ref{eq:gen_opt}), we employ a gradient-based approach
that incorporates sensitivity information to accelerate convergence to an optimal solution. Specifically, the sample average approximation (SAA) method, also referred to as sample path optimization, is used to optimize the mean and covariance of the random field\cite{rubinstein_shapiro}. We limit our attention to the special case where each objective and constraint functions are equal to the mean of an output functional, since this special case encompasses the problems of interest in this work.

In the SAA method, the objective functions and constraints are approximated using the Monte Carlo method. For example, the mean of the functional $F(e)$ is estimated as
\begin{equation}
\mathbb{E}[F] \approx \frac{1}{N} \sum_{n=1}^{N} F_n.
\label{eq:mc}
\end{equation}
The process for propagating distributed uncertainty to the quantities of interest is summarized below:
\begin{enumerate}
\item Generate a $N \times N_{KL}$ matrix of independent Gaussian random variables.
\item For each Monte Carlo sample, construct a realization of the random field $e_n(x)$ using the K-L expansion (\ref{eq:kl}).
\item Evaluate the functional of interest $F_n = F(e_n)$ for each realization.
\item Estimate the moments of $F$ according to Equation (\ref{eq:mc}).
\end{enumerate}
The convergence rate of the Monte Carlo estimate (\ref{eq:mc}) is $O(N^{-1/2})$, and therefore a large number of Monte Carlo samples are typically required. However, the Monte Carlo samples can be evaluated in parallel, greatly reducing the time required to evaluate Equation (\ref{eq:mc}).

The SAA method transforms the stochastic optimization problem (\ref{eq:gen_opt}) into a deterministic optimization problem. This is achieved by fixing the set of realizations $\{\bxi_n\}_{n=1}^{N}$ of the random input vector used to compute the Monte Carlo estimates of the objective and constraint functions. The SAA method therefore solves the following modified optimization problem, where the objective and constraint functions have been replaced by their Monte Carlo estimates:
\begin{equation}
\begin{aligned}
& \hat{\bp}_N^* =
& & \underset{\bp \in \mathcal{P}}{\arg \min}
& & \hat{f}_N(\bp) \\
& & & \ \ \ \ \text{s.t.}
& &  \hat{g}_N(\bp) \leq 0 \\
& & & & & \hat{h}_N(\bp) = 0.
\end{aligned}
\label{eq:approx_opt}
\end{equation}
The subscript $N$ has been added to emphasize the number of samples used to construct the estimators. The deterministic optimization problem that results from fixing the samples can be solved iteratively to update the solution, using the same set of realizations $\{\bxi_n\}_{n=1}^{N}$ at each iteration. The solution of the deterministic optimization problem, denoted $\hat{\bp}^*_N$, is an estimator of the true solution $\bp^*$. 

In the unconstrained case, $\hat{f}_N(\hat{\bp}^*_N) \rightarrow f(\bp^*)$ and $\hat{\bp}^*_N \rightarrow \bp^*$ as $N\rightarrow \infty$ with probability one if $\bp^*$ is a unique minimizer of $f$ and the family $\{|F(\samp,\bp)|,\,\bp \in \mathcal{P}\}$ is dominated by a measurable function, i.e. if there exists a measurable function $G(\samp)$ such that  $|F(\samp,\bp)| \leq G(\samp)$ for all points $\samp \in \Samp$\cite{rubinstein_shapiro}. Moreover, if the families $\{||\nabla F(\samp,\bp)||,\,\bp \in \mathcal{P}\}$ and $\{||\nabla^2 F(\samp,\bp)||,\,\bp \in \mathcal{P}\}$ are dominated by measurable functions, then, assuming the Hessian matrix $\mathbf{B} = \mathbb{E}[\nabla^2 F(\samp,\bp^*)]$ is nonsingular, 
\begin{equation}
N^{1/2}(\hat{\bp}^*_N - \bp^*) \overset{\text{i.d.}}{\rightarrow} \mathcal{N}(0,\mathbf{B}^{-1}\boldsymbol\Sigma\mathbf{B}^{-1}),
\end{equation}
\begin{equation}
N^{1/2}(\hat{f}_N(\hat{\bp}^*_N) - f(\bp^*)) \overset{\text{i.d.}}{\rightarrow} \mathcal{N}(0,\gamma^2),
\end{equation}
where $\overset{\text{i.d.}}{\rightarrow}$ represents convergence in distribution and 
\begin{equation}
\boldsymbol\Sigma = \mathbb{E}[\nabla F(\samp,\bp^*) \nabla F(\samp,\bp^*)^\intercal],
\end{equation}
\begin{equation}
\gamma^2 = \mathbb{E}[F(\samp,\bp^*)^2] - f(\bp^*)^2.
\end{equation}
Thus, the SAA approximate solution and approximate objective function converge like $N^{-1/2}$. Since the true solution $\bp^*$ is unknown, the quantities $\mathbf{B}$, $\boldsymbol\Sigma$ and $\gamma^2$ are replaced by consistent estimates computed from the same realizations $\{\bxi_n\}_{n=1}^{N}$ used to solve the problem\cite{rubinstein_melamed}. It is also possible to assess the SAA solution quality by constructing a confidence bound on the optimality gap $f(\hat{\bp}_N^*) - f(\bp^*)$\cite{mak_1999}.

The reduction of the stochastic optimization problem into a deterministic optimization problem allows for the use of one of many algorithms designed for the efficient solution of deterministic optimization problems. Thus, the SAA method is well-suited to solving constrained stochastic optimization problems. A convergence rate of $N^{-1/2}$ for the constrained problem can also be observed under certain conditions\cite{rubinstein_shapiro}. Numerous methods have been devised for solving deterministic optimization problems with both nonlinear objectives and nonlinear constraints. One such method, the sequential quadratic programming (SQP) method, is reviewed next.  

\subsection{Sequential quadratic programming}

An efficient approach to solving (\ref{eq:approx_opt}) is the sequential quadratic programming method. Given an approximate solution $\hat{\bp}^k$, the SQP solves a quadratic programming subproblem to obtain an improved approximate solution $\hat{\bp}^{k+1}$. This process is repeated to construct a sequence of approximations that converge to a solution $\hat{\bp}^*$\cite{boggs_1996}. The quadratic subproblems are formed by first constructing the Lagrangian function from the objective and constraint functions. A quadratic objective is constructed from the second-order Taylor series expansion of the Lagrangian, and the constraints are replaced with their linearizations. The solution of the quadratic subproblem produces a search direction, and a linesearch can be applied to update the approximate solution. To construct the second-order Taylor series of the Lagrangian, the Hessian is estimated using a quasi-Newton update formula, such as the Broyden-Fletcher-Goldfarb-Shanno (BFGS) formula\cite{nocedal}. Local convergence of the SQP algorithm  requires that the initial approximate solution is close to a local optimum and that the approximate Hessian is close to the true Hessian. Global convergence requires sufficient decrease in a merit function that measures the progress towards an optimum. More details on local and global convergence of SQP methods can be found in \cite{boggs_1996}.

\section{Sensitivity analysis of Gaussian random fields}

In this section, we perform sensitivity analysis of a system's output statistics with respect to the mean and covariance of Gaussian random field input uncertainty. This sensitivity information is used to optimize the mean and covariance functions via the SAA method described in the previous section.

\subsection{Pathwise sensitivities}

To compute the sensitivity of an output statistic, e.g. $\nabla_p \E[F(p)]$, we use the pathwise sensitivity method. The pathwise sensitivity method relies upon interchanging the differentiation and expectation operators. For example, to compute an unbiased estimator of the gradient of $f = \mathbb{E}[F(\samp,\bp)]$ with respect to a parameter $p$, we simply interchange differentiation and integration:
\begin{equation}
\frac{\partial}{\partial p}\mathbb{E}[F(\samp,\bp)] = \mathbb{E}\left[\frac{\partial}{\partial p}F(\samp,\bp)\right].
\label{eq:pathwise}
\end{equation}
Sufficient conditions that allow for this interchange will be discussed subsequently. 

The pathwise sensitivity method can applied directly to the Monte Carlo estimate of $\mathbb{E}[F(\samp,\bp)]$. Replacing the expectation with its Monte Carlo estimate, and exchanging summation and differentiation gives
\begin{equation}
\mathbb{E}\left[\frac{\partial F}{\partial p} \right] \approx \frac{\partial \hat{f}_N}{\partial p} = \frac{1}{N}\sum_{n=1}^{N}\frac{\partial F_n}{\partial p}.
\label{eq:pathwise_mc}
\end{equation}
In the context of the SAA method, the derivatives $\partial F_n/\partial p$ represent the sensitivity of the random functional $F(\samp,\bp)$ for a particular realization of the random field $e_n \equiv e(x,\boldsymbol\xi_n)$ where all random inputs are held fixed. To compute the sensitivity $\partial F_n/\partial p$, we first apply the chain rule to rewrite this sensitivity:
\begin{equation}
\frac{\partial F_n}{\partial p} = \frac{\partial F_n}{\partial e_n}\frac{\partial e_n}{\partial p}.
\end{equation}
If the functional $F$ depends explicitly on the random field $e$, the derivative $\partial F_n/\partial e_n$ can be computed directly. As mentioned previously, $F$ may alternatively be a functional of the solution of some system of equations depending on $e$. In that case, the derivative $\partial F_n/\partial e_n$ can be computed efficiently using the adjoint method\cite{giles_2000}. We now turn our attention to computing the sample path sensitivity $\partial e_n/\partial p$.

\subsection{Sample path sensitivities}

We consider computing the sensitivity of the sample path $e_n \equiv e(x,\boldsymbol\xi_n; \bp_m, \bp_c)$ with respect to the parameters which control the mean and covariance of the random field, i.e. the $\bp_m$ and $\bp_c$ introduced previously. The sensitivity of the sample path with respect to any parameter $p_m$ controlling the mean can be analytically derived from the K-L expansion given by Equation (\ref{eq:kl}). Since the eigenvalues and eigenvectors in the K-L expansion are independent of $p_m$, only the first term in the K-L expansion depends on $p_m$. Thus, we have 
\begin{equation}
\frac{\partial e_n}{\partial p_m} = \frac{\partial \bar{e}}{\partial p_m}.
\end{equation}

Computing the sensitivity of the sample path with respect to a parameter $p_c$ controlling the covariance is more involved. The pathwise sensitivity method has typically been applied to problems in computational finance and chemical kinetics where the sample paths of the random process can be differentiated analytically with respect to the parameters of interest\cite{broadie_1996,sheppard}. However, the sample path sensitivity of a random field can not, in general, be differentiated analytically with respect a parameter controlling the covariance matrix. For a Gaussian random field, we can use its K-L expansion to compute these sensitivities using eigenvalue/eigenvector perturbation theory. We focus on computing the sensitivities of the discretized random field, since numerical computation of the pathwise sensitivity estimate is the ultimate goal. We first consider the general case of computing the sensitivity of the sample path with respect to a covariance parameter $p_c$, and then the special case where the parameter of interest controls the variance of a random field with fixed correlation function.

\subsubsection{General case}

Since the covariance matrix is a function of $p_c$, its eigenvalues and eigenvectors are also functions of $p_c$. Applying the chain rule to the Karhunen-Lo\`{e}ve expansion, we have
\begin{equation}
\frac{\partial e_n}{\partial p_c} = \sum_{i=1}^{N_{KL}}\left( \frac{1}{2\sqrt{\lambda_i}} \phi_i\frac{\partial \lambda_i}{\partial p_c} + \sqrt{\lambda_i}\frac{\partial \phi_i}{\partial p_c} \right)\xi_i(\samp_n).
\label{eq:sample_deriv}
\end{equation}
Note that since the pathwise sensitivity approach is used, the random variables $\xi_i(\samp_n)$ remain fixed. Equation (\ref{eq:sample_deriv}) is only valid if the eigenvalues and eigenvectors in the K-L expansion are differentiable functions of $p_c$. It can be shown, via the implicit function theorem, that if the eigenvalues of $\mathbf{C}$ are simple (i.e., have algebraic multiplicity one), then the eigenvalues and eigenvectors of $\mathbf{C}$ are infinitely differentiable with respect to $p_c$\cite{magnus_2007}. % See Theorem 7 p. 179 of magnus_2007. 
If the eigenvalues remain simple as $p_c$ is varied over some range of values, then the eigenvalues and eigenvectors are differentiable over that range of $p_c$.

For an arbitrarily chosen covariance matrix, varying a parameter $p_c$ controlling the covariance function is unlikely the result in duplicate eigenvalues. To see this, first note that the difference between the dimension of the space of $n \times n$ symmetric positive definite matrices and the dimension of the subspace of $n \times n$ symmetric positive definite matrices with repeated eigenvalues is at least two, which can be proved using a simple counting argument\cite{lax}. A curve in $N$ dimensional space is unlikely to pass through a $N-2$ dimensional subspace, e.g. an arbitrary curve in the plane ($N = 2$) is unlikely to pass through a given point in that plane. This gives rise to the ``avoidance of crossing'' phenomena: as $p_c$ is varied, the eigenvalues of a symmetric matrix are extremely unlikely to cross, and thus are likely to remain simple\cite{lax}. Thus, the eigenvalues and eigenvectors are likely to remain differentiable functions of $p_c$ as $p_c$ is varied. 

\begin{figure}[htbp]
\centering
\vspace{0.2cm}
\includegraphics[width=0.6\textwidth]{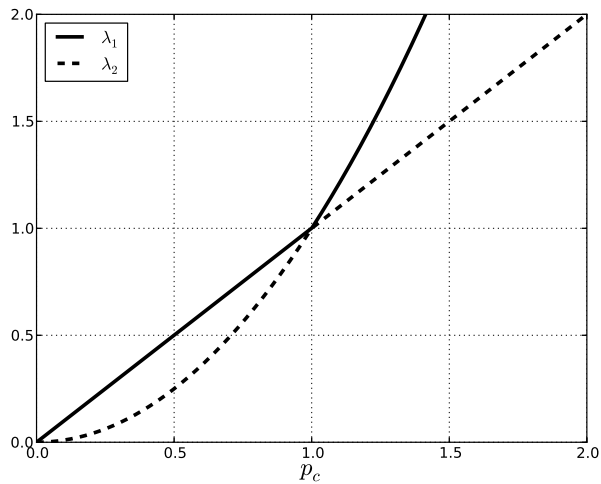}
\caption{Example of crossing eigenvalues of a symmetric, positive definite matrix.}
\label{fig:eig_cross}
\end{figure}

Of course, it is easy to design cases where the eigenvalues cross. For example, consider the matrix 
\begin{equation}
\mathbf{C} = \left[\begin{array}{cc} 
p_c & 0 \\
0 & p_c^2
\end{array}\right]
\end{equation}
over the range $p_c \in (0,\infty)$. The eigenvalues of this matrix are plotted in Figure \ref{fig:eig_cross}, which clearly shows the two eigenvalues crossing at $p_c = 1$. At the point of crossing, the eigenvalues are not differentiable with respect to $p_c$, which can be visualized by the ``kinks'' in the two curves at $p_c = 1$. However, such cases are extremely unlikely to occur for arbitrary covariance matrices, where the elements are not deliberately chosen to produce crossing eigenvalues.

When the eigenvalues are simple, the derivatives of the eigenvalues and eigenvectors can be computed using established results from eigenvalue perturbation theory:
\begin{equation}
\frac{\partial \lambda_i}{\partial p_c} = \phi_i^\intercal\frac{\partial \mathbf{C}}{\partial p_c}\phi_i,
\end{equation}
and
\begin{equation}
\frac{\partial \phi_i}{\partial p_c} = -(\mathbf{C} - \lambda_i \mathbf{I})^+\frac{\partial \mathbf{C}}{\partial p_c}\phi_i,
\end{equation}
where $(\mathbf{C} - \lambda_i \mathbf{I})^+$ denotes the Moore-Penrose pseudoinverse of the matrix $(\mathbf{C} - \lambda_i \mathbf{I})$\cite{deleeuw}. Since the explicit dependence of the entries of the covariance matrix $\mathbf{C}$ on $p_c$ is assumed to be known, the sensitivities of the eigenvalues and discretized eigenvectors in the K-L expansion can be computed in closed form.

One practical issue that arises when using the pathwise sensitivity method results from the sign ambiguity of the eigenvectors. Specifically, although the eigenvector $\phi_i(p_c)$ is differentiable with respect to $p_c$ (and therefore continuous), perturbing $p_c$ by some small $\varepsilon$ may result in $\phi_i(p_c + \varepsilon)$ being very different from $\phi_i(p_c)$ as a result of sign ambiguity. This issue is resolved by choosing the sign that results in the ``closer'' eigenvector: if $\|\phi_i(p_c + \varepsilon) + \phi_i(p_c)\|_2 < \|\phi_i(p_c + \varepsilon) - \phi_i(p_c)\|_2$, then the sign of $\phi_i(p_c + \varepsilon)$ is flipped.

\subsubsection{Special case: sensitivity with respect to the variance}

Computing the sample path sensitivities can be simplified if the parameter $p_c$ only scales the variance of the random field, but does not change its correlation function. Consider a random field $\tilde{e}(x,\samp)$ with unit variance, i.e. $\mathbb{E}[\tilde{e}^2(x,\samp)] = 1$ everywhere. Scaling this random field by the function $\sigma(x)$ produces the random field $e(x,\samp) = \sigma(x)\tilde{e}(x,\samp)$ with non-stationary variance $\sigma^2(x)$\cite{abrahamsen}. The covariance function of the process $\tilde{e}(x,\samp)$, denoted $\rho(x_1,x_2)$, satisfies the property $x_1 = x_2 \implies \rho(x_1,x_2) = 1$. The corresponding covariance function of the scaled process $e(x,\samp)$ is given by $\mathbf{C}(x_1,x_2) = \sigma(x_1)\sigma(x_2)\rho(x_1,x_2)$. 

Suppose the function $\sigma(x)$ depends smoothly on the parameters $\bp_c$. Rather than simulating the random field $e(x,\samp)$ with non-stationary variance, we instead simulate the unit variance field $\tilde{e}(x,\samp)$ and set $e_n(x) = \sigma(x)\tilde{e}_n(x)$. Then, the sample path sensitivity with respect to $p_c$ can be computed as
\begin{equation}
\frac{\partial e_n}{\partial p_c} = \frac{\partial e_n}{\partial \sigma}\frac{\partial \sigma}{\partial p_c} = \tilde{e}_n\frac{\partial \sigma}{\partial p_c}.
\end{equation}
This greatly simplifies the sensitivity calculation since the K-L expansion only needs to be computed once. This eliminates the issues caused by the sign ambiguity of the eigenvectors since the same set of eigenvectors are used throughout the optimization. The computational cost of performing optimization with this approach is also lower since it does not require the sensitivity of the K-L expansion to be computed at each optimization step. However, this difference in computational cost may be small compared to the cost of computing the objective and constraint function estimates, which typically require many Monte Carlo simulations to be performed. If each Monte Carlo sample is computationally expensive, e.g. requires solving a system of partial differential equations, then the relative savings will be very small.

\begin{figure}[htbp]
\centering
\subfloat[Stationary standard deviation (left) and resulting realizations (right)]{\includegraphics[width=0.49\textwidth]{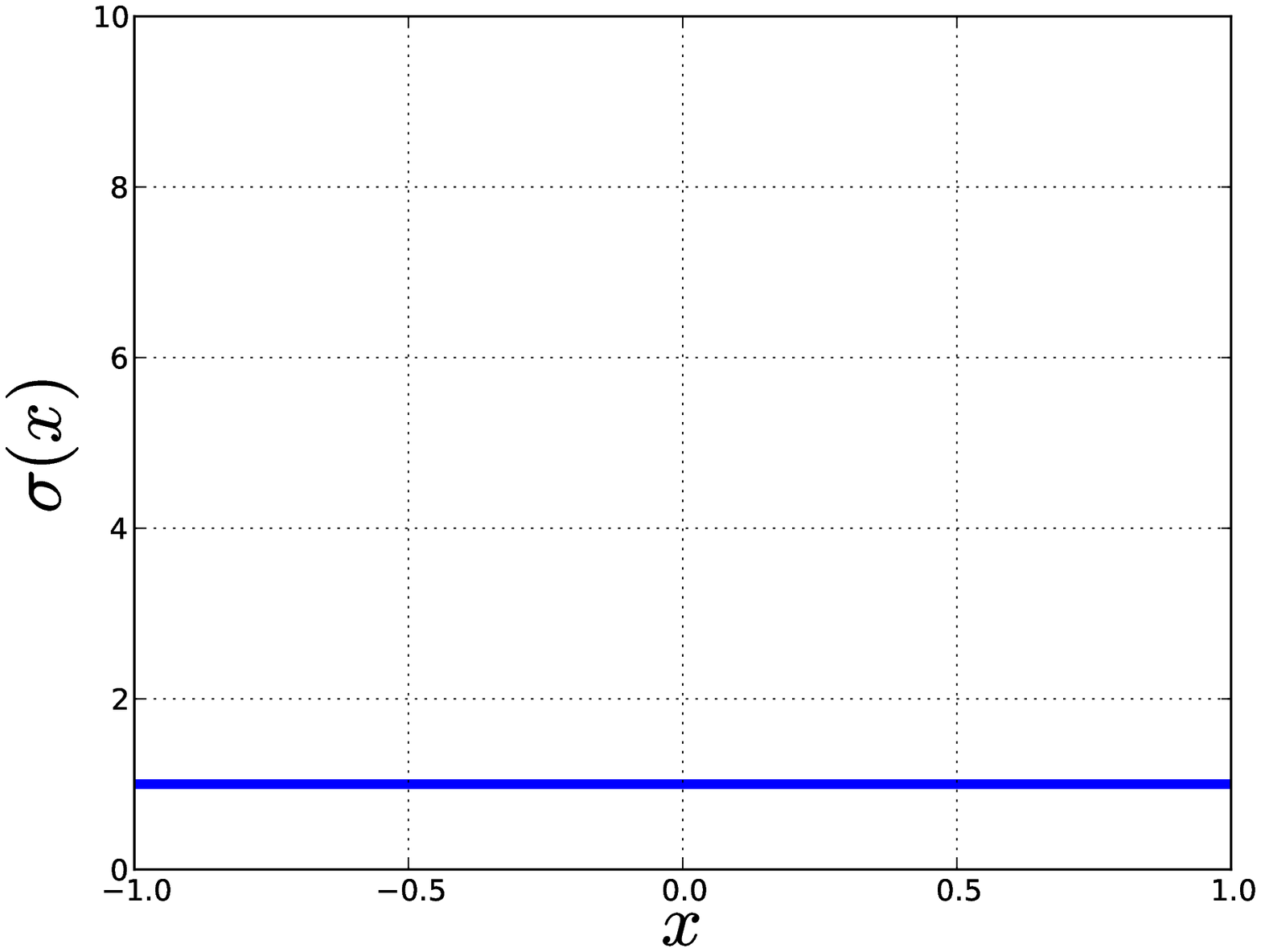}
\includegraphics[width=0.49\textwidth]{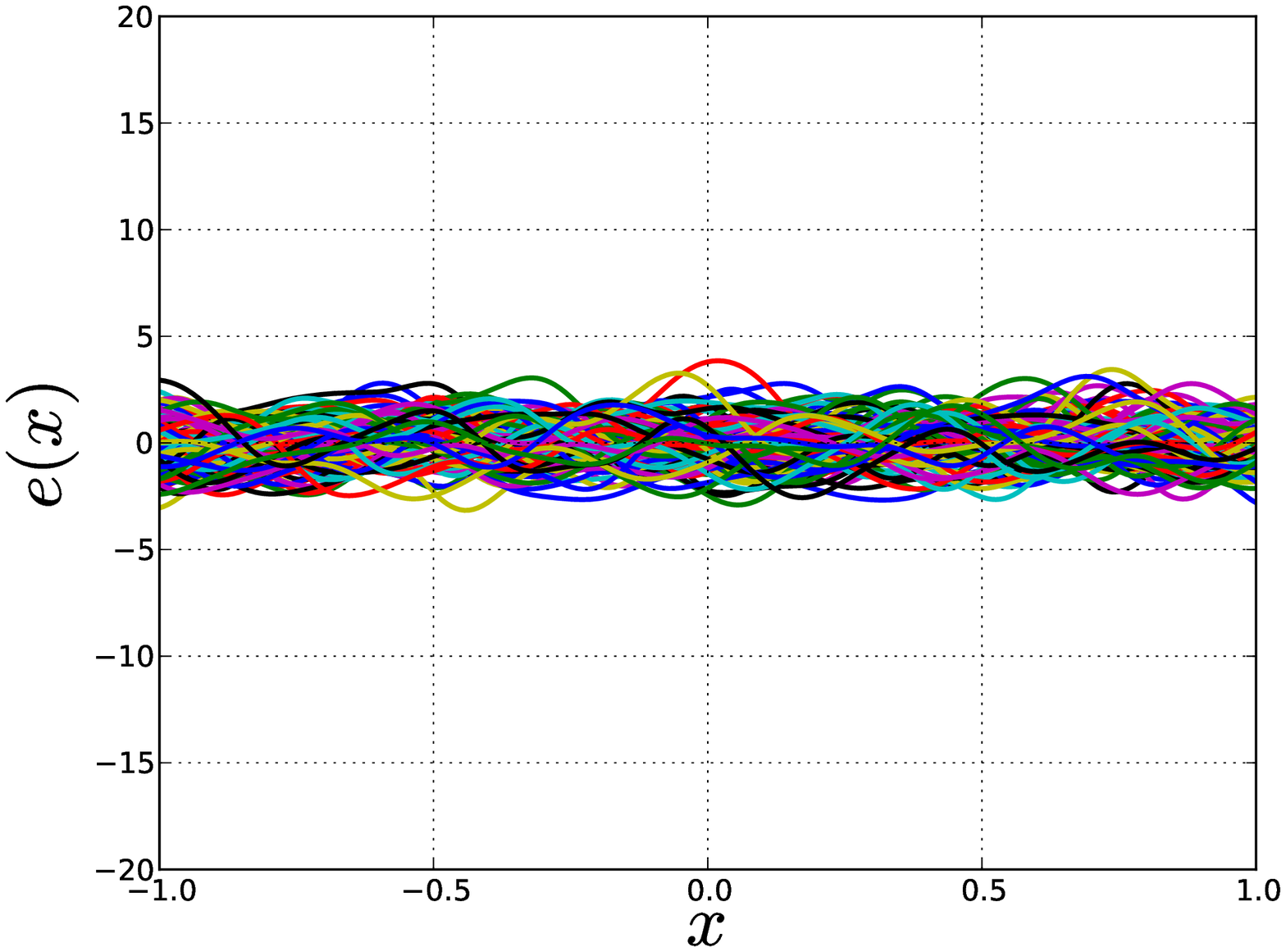}}
\newline
\subfloat[Non-stationary standard deviation (left) and resulting realizations (right)]{\includegraphics[width=0.49\textwidth]{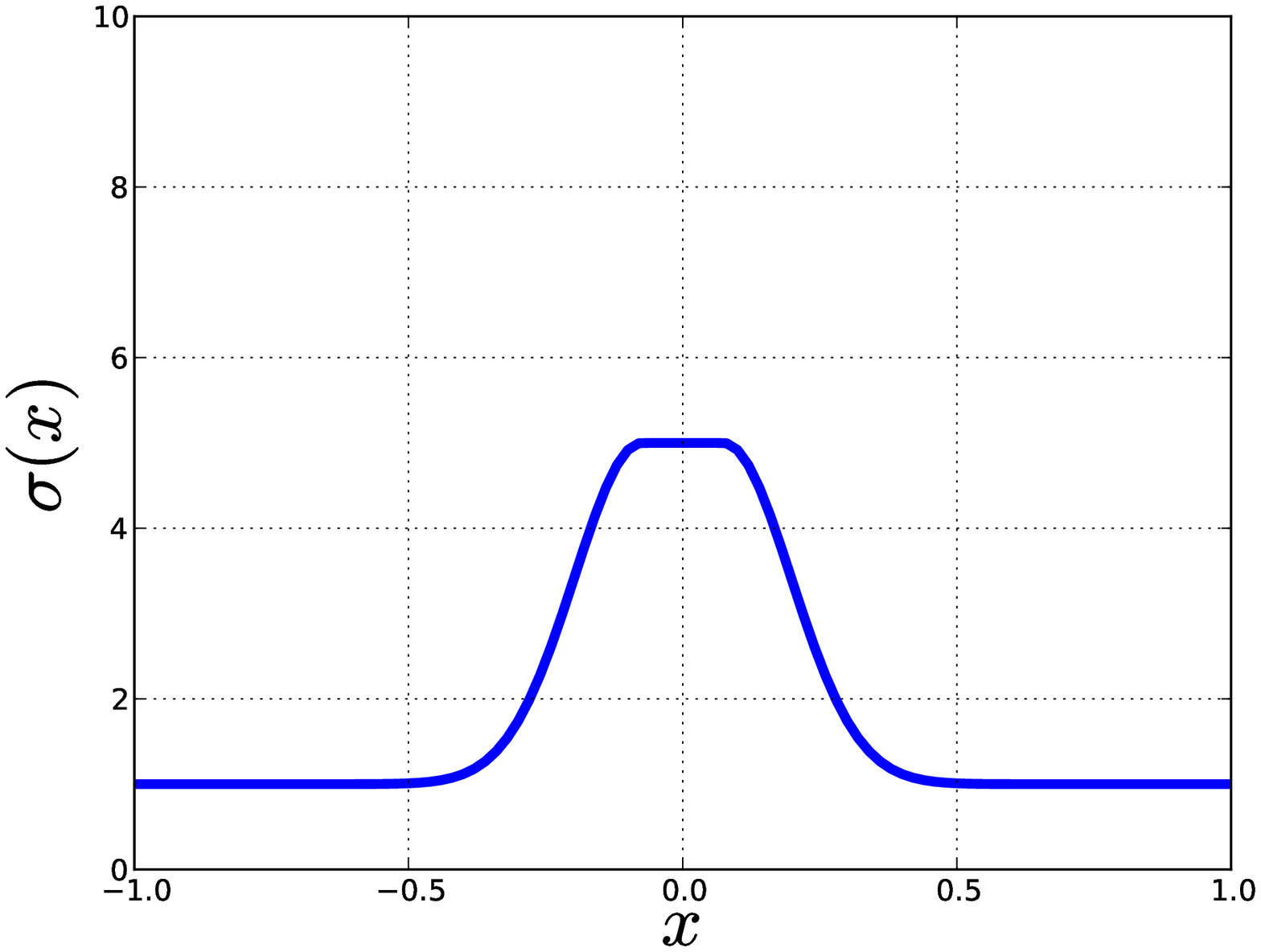}
\includegraphics[width=0.49\textwidth]{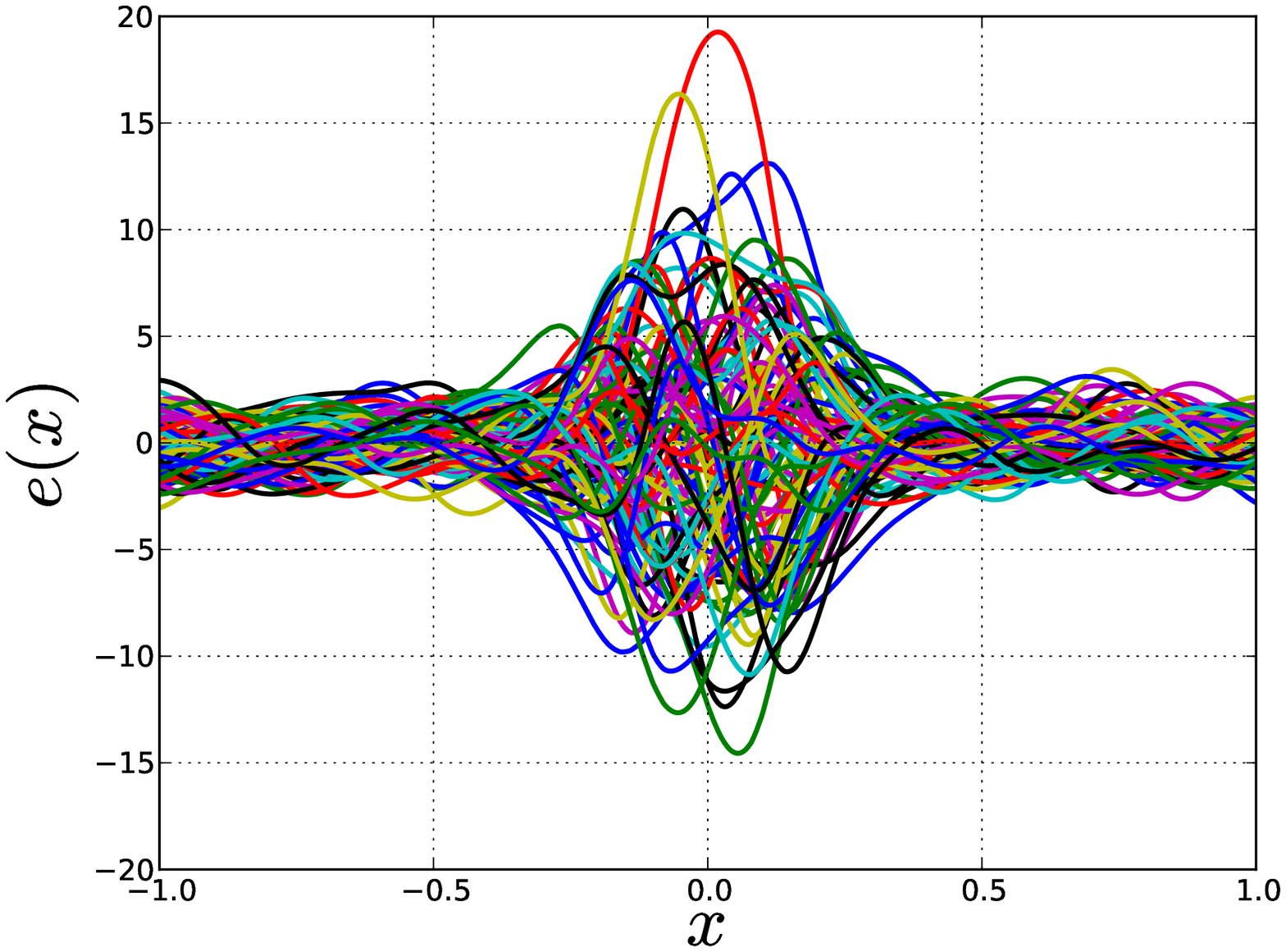}}
\caption{Random fields with stationary (top) and non-stationary (bottom) standard deviation.}
\label{fig:rf_nonstat_var}
\end{figure}

Figure \ref{fig:rf_nonstat_var} illustrates scaling a random field with stationary variance to produce realizations of a random field with a spatially varying variance. The original random field, shown at the top, is a zero-mean Gaussian random field with a squared exponential covariance function. The scaled random field, shown on the bottom, is also a zero-mean Gaussian random field. However, the increase in the standard deviation near $x = 0$ produces realizations with more variability in this region than the original random field with stationary variance.

\subsection{Interchanging differentiation and expectation}

As mentioned previously, applying the pathwise sensitivity method requires that the interchange of differentiation and integration is justified. We now address which conditions on $F$ and $\bp$ ensure that this interchange is justified. The first requirement is that the random vector $\boldsymbol\xi$ must be independent of the parameters $\bp$. Since we use the K-L expansion to simulate the random field, this is true by construction: changing the parameters $\bp$ only changes the eigenvalues and eigenvectors in the K-L expansion, thus the random vector $\boldsymbol\xi$ is independent of the parameters $\bp$.

The second requirement is on the regularity of the function $F(\samp,p)$ (for simplicity, we only one parameter $p$). Interchanging differentiation and integration requires that the following interchange of limit and integration is justified:
\begin{equation}
\mathbb{E}\left[\lim_{h\rightarrow 0}\frac{F(\samp,p + h) - F(\samp,p)}{h}\right] = \lim_{h\rightarrow 0} \mathbb{E}\left[\frac{F(\samp,p + h) - F(\samp,p)}{h}\right].
\label{eq:lim_exchange}
\end{equation}
A necessary and sufficient condition for this interchange to be valid is that the difference quotients $Q_h = h^{-1}[F(\samp,p+h) - F(\samp,p)]$ are uniformly integrable, i.e. that
\begin{equation}
\lim_{c \rightarrow \infty} \sup_{h} \mathbb{E} [|Q_h| \mathbf{1}\{|Q_h| > c\}] = 0,
\end{equation}
where $\mathbf{1}\{|Q_h| > c\}$ is the indicator function. This condition is not readily verified for practical problems, since the analytical distribution of $F$ is typically unknown. We instead provide a set of sufficient conditions that are more straightforward to verify in practice, following reference \cite{glasserman}. Recall that $F$ is a functional of the random field $e(\samp,x;p)$, and denote by $D_F \subset \mathbb{R}^{|\Samp|}$ the set of points in $\Samp$ where $F$ is differentiable with respect to $e$. The following are sufficient conditions for the interchange of the limit and expectation in (\ref{eq:lim_exchange}).
\begin{itemize}[label={},leftmargin=*]
\item{\textbf{(A1)}} For every $p \in \mathcal{P}$ and $x \in X$, $\partial e(x,\samp;p)/\partial p$ exists with probability 1.
\item{\textbf{(A2)}} For every $p \in \mathcal{P}$, $\mathbb{P}[e(x,\samp;p) \in D_F] = 1$.
\item{\textbf{(A3)}} $F$ is Lipschitz continuous, i.e. there exists a constant $k_F < \infty$ such that for all $u(x)$, $v(x)$,
\begin{equation}
|F(u) - F(v)| \leq k_F\|u - v\|.
\end{equation}
\item{\textbf{(A4)}} For every $x \in X$, there exists a random variable $k_e$ such that for all $p_1,\, p_2\in \mathcal{P}$,
\begin{equation}
|e(x,\samp;p_2) - e(x,\samp;p_1)| \leq k_e|p_2 - p_1|,
\end{equation}
and $\mathbb{E}[k_e] < \infty$.
\end{itemize}
Conditions (A3) and (A4) imply that $F$ is Lipschitz continuous in $p$ with probability one. Taking $\kappa_F = k_F \sup_x k_e$,
\begin{equation}
|F(\samp,p_2) - F(\samp,p_1)| \leq \kappa_F |p_2 - p_1|.
\end{equation}
We can then bound the difference quotient:
\begin{equation}
\left|\frac{F(\samp,p + h) - F(\samp,p)}{h} \right| \leq \kappa_F,
\end{equation}
and apply the dominated convergence theorem to interchange the expectation and limit in (\ref{eq:lim_exchange}). Thus, conditions (A1)-(A4) are sufficient conditions for the pathwise sensitivity estimate to be unbiased.

Conditions (A3) and (A4) together determine if $F$ is almost surely Lipschitz continuous, and thus determine what type of input parameters and output quantities of interest can be treated with the pathwise sensitivity method. The previous section gave conditions for the differentiability of the sample paths, i.e. that the covariance function depends smoothly on $\bp$ and have simple eigenvalues. Output functionals that may change discontinuously when smooth perturbations are made to the random field are not Lipschitz continuous almost surely. Thus, condition (A3) excludes failure probabilities, e.g. $\mathbb{P}(F \geq c) = \mathbb{E}[\mathbf{1}\{F \geq c\}]$, since the indicator function $\mathbf{1}\{F \geq c\}$ is discontinuous when $F = c$. This difficulty can be remedied to some degree using a smoothed version of the indicator function, but this introduces additional error to the sensitivity\cite{foque_2003}. Conditions (A2) and (A3) do permit functions that fail to be differentiable at certain points, as long as the points at which differentiability fails occur with probability zero, and $F$ is continuous at these points.

\section{Application: variance optimization}

To demonstrate the proposed optimization framework, we consider an optimization problem with the design variables being the variance of a random field. The random field $e(x,\samp)$ is defined on the domain $X = [0,1]$ and has a squared exponential correlation function:
\begin{equation}
\rho(x_1,x_2) = \exp\left[-\frac{(x_1 - x_2)^2}{2L^2}\right],
\end{equation}
with correlation length $L = 0.1$. The standard deviation $\sigma(x)$ of the random field is a spatially dependent function. We seek to minimize the sum of two competing cost functions that depend on $\sigma(x)$ as a (spatially varying) parameter. The first cost function penalizes variability:
\begin{equation}
f_1 = \mathbb{E}\left[ \int_0^1 e^2(x,\samp)w(x)\ dx \right],
\label{eq:f1}
\end{equation}
where $w(x)$ is a non-negative weighting function. The weighting function  specifies which regions are most sensitive to increased variability. Regions where $w(x)$ is large correspond to regions where variability has the largest impact on the system. The second cost function is inversely proportional to the variability:
\begin{equation}
f_2 = \int_0^1 \frac{1}{\sigma(x)}\ dx
\label{eq:f2}
\end{equation}
We seek to determine the standard deviation field $\sigma^*(x)$ that minimizes the sum of the two cost functions:
\begin{equation}
\begin{aligned}
\sigma^*(x) & = \underset{\sigma(x)}{\arg\min} &  f_1 + f_2 \equiv f % \\
% &\ \ \ \ \ \ \  \text{s.t.} & \sigma(x) > 0.
\end{aligned}
\label{eq:opt_state}
\end{equation}
This model problem is analogous to a tolerance optimization problem. Reducing tolerances (thereby increasing the variance $\sigma^2(x)$) can improve the performance of the system. This behavior is reflected in the cost function $f_1$. Moreover, certain regions of the domain are more sensitive to variability than others, as expressed by the weight function $w(x)$. On the other hand, it is costly to reduce tolerances, and the cost of reducing tolerances increases monotonically, as reflected in the form of $f_2$.

The optimal solution to (\ref{eq:opt_state}) can be derived analytically using the calculus of variations. The expectation and spatial integration can be interchanged in Equation (\ref{eq:f1}) to give
\begin{equation}
f_1 = \int_0^1 \mathbb{E}[e^2(x,\samp)]w(x)\ dx = \int_0^1 \sigma^2(x)w(x)\ dx.
\label{eq:f1_expand}
\end{equation}
The first variation of $f$ can then be computed directly:
\begin{equation}
\delta f = \int_0^1 \left(2\sigma(x) w(x) - \frac{1}{\sigma^2(x)}\right)\delta \sigma(x)\ dx.
\end{equation}
Enforcing stationarity by setting $\delta f = 0$, the optimal standard deviation field is found to be
\begin{equation}
\sigma^*(x) = \left[\frac{1}{2w(x)} \right]^{1/3}.
\end{equation}
Note that this optimal is unique since both $f_1$ and $f_2$ are strictly convex functionals.

As an example, we choose the weight function to be $w(x) = 2 + \sin(2\pi x)$. The standard deviation field is discretized with $N_\sigma = 20$ cubic B-spline basis functions $B_i$:
\begin{equation}
\sigma(x) = \sum_{i = 1}^{N_\sigma}\sigma_i B_{i}(x).
\end{equation}
To demonstrate our method, the Monte Carlo method is used to compute an unbiased estimate of $f_1$, rather than computing it directly from Equation (\ref{eq:f1_expand}):
\begin{equation}
\hat{f}_1 = \frac{1}{N} \sum_{n=1}^{N} \int_0^1 e^2_n(x) w(x) \ dx
\end{equation}
For each Monte Carlo sample, the integral is evaluated using composite Gaussian quadrature with 20 intervals and a third order rule on each interval. The same quadrature rule is used to compute $f_2$. The SAA equivalent of (\ref{eq:opt_state}) results from replacing the objective function $f_1$ with its unbiased estimate:
\begin{equation*}
\begin{aligned}
\hat{\sigma}^*(x) & = \underset{\sigma(x)}{\arg\min} &  \hat{f}_1 + f_2 % \\
% &\ \ \ \ \ \ \  \text{s.t.} & \sigma(x) > 0.
\end{aligned}
\end{equation*}
This optimization problem is solved using the SQP algorithm with a BFGS update to approximate the Hessian as implemented in the NLopt package\cite{nlopt}. The pathwise estimate of the sensitivity $\partial f_1/\partial \sigma(x)$, which is an unbiased estimate of the true gradient, is computed as 
\begin{equation}
\frac{\partial \hat{f}_1}{\partial \sigma} = \frac{1}{N} \sum_{n=1}^{N} \int_0^1 2 w(x) e_n(x) \frac{\partial e_n}{\partial \sigma}\ dx
\end{equation}
The sample path sensitivity $\partial e_n/\partial \sigma$ can be computed using either approach described previously, i.e. by computing the sensitivity of the K-L expansion or by computing sensitivities for a unit-variance random field scaled by $\sigma(x)$. We use both approaches to compare their effectiveness.

\begin{figure}[htbp]
\centering
\subfloat[$N = 10^2$]{\includegraphics[width=0.49\textwidth]{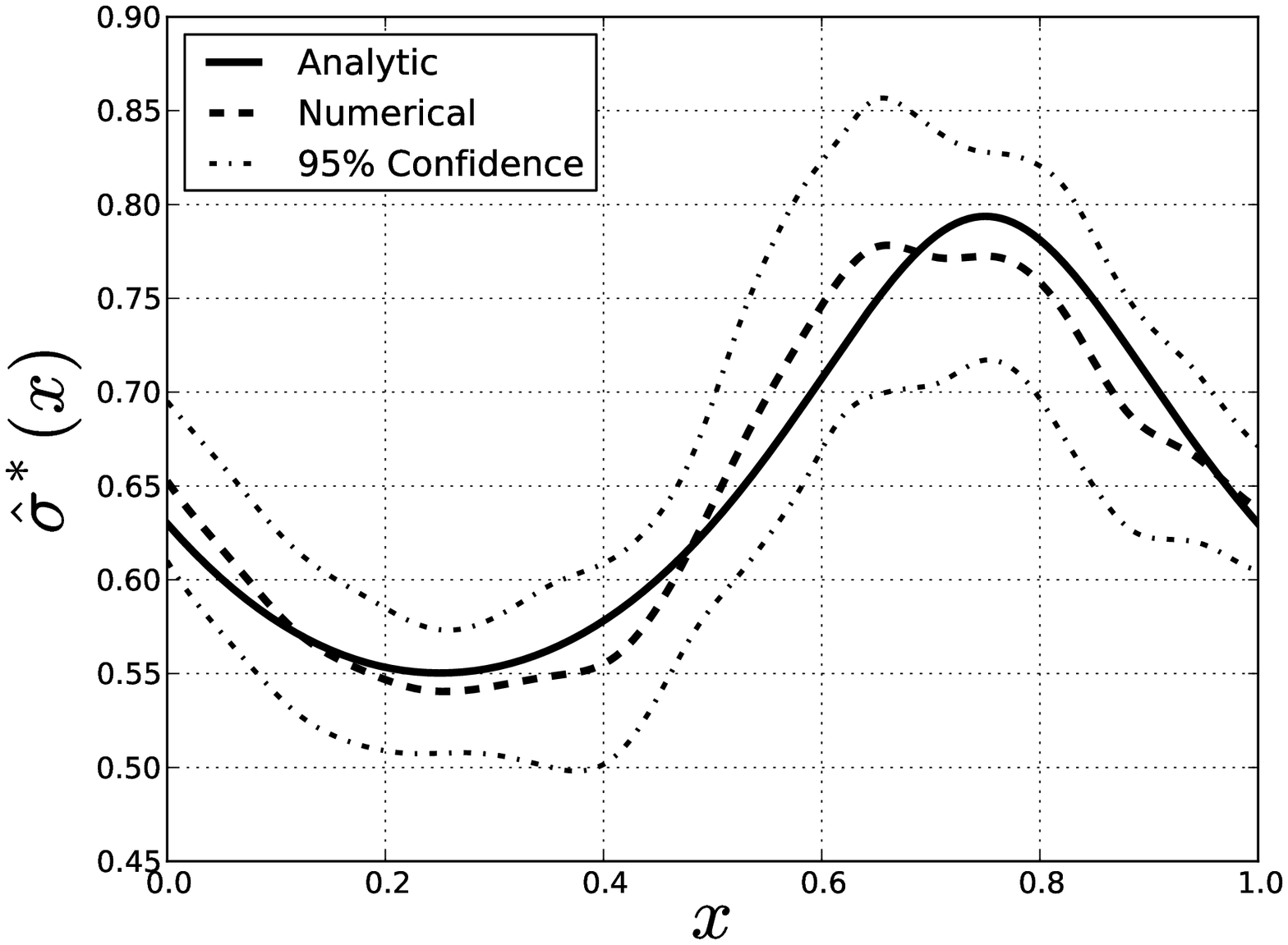}}
\hspace{0.1cm}
\subfloat[$N = 10^3$]{\includegraphics[width=0.49\textwidth]{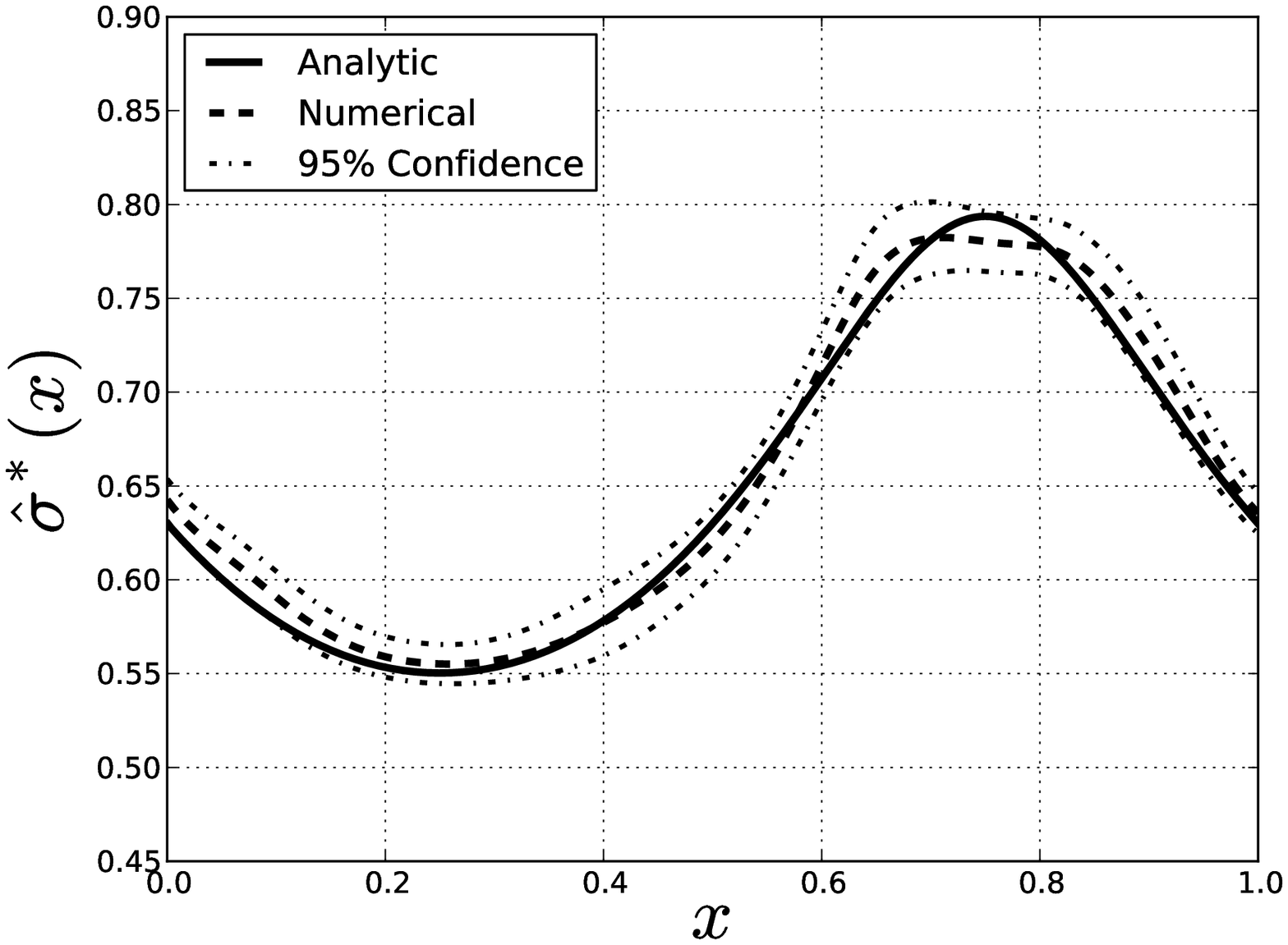}}
\newline
\subfloat[$N = 10^4$]{\includegraphics[width=0.49\textwidth]{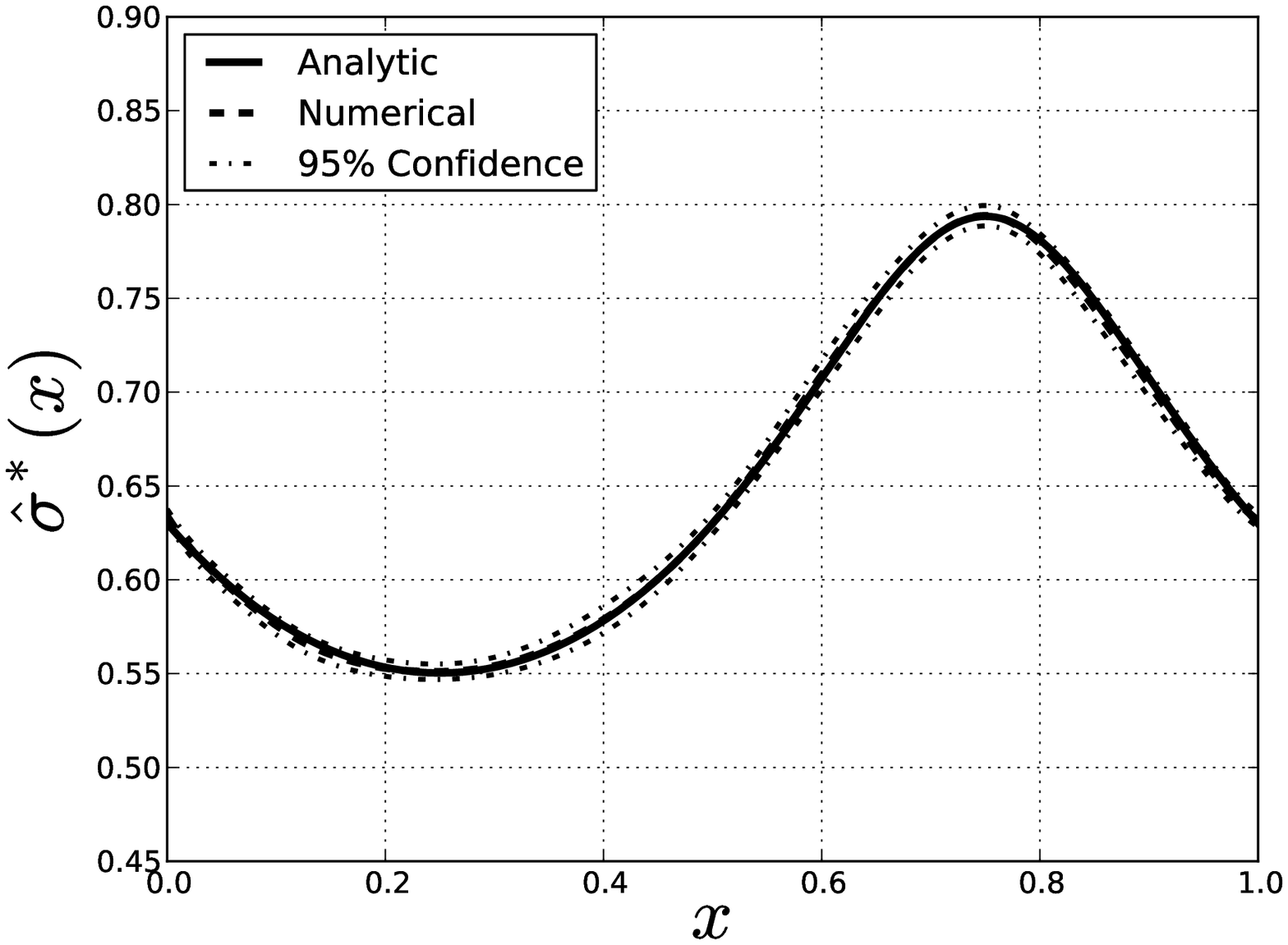}}
\hspace{0.1cm}
\subfloat[$N = 10^5$]{\includegraphics[width=0.49\textwidth]{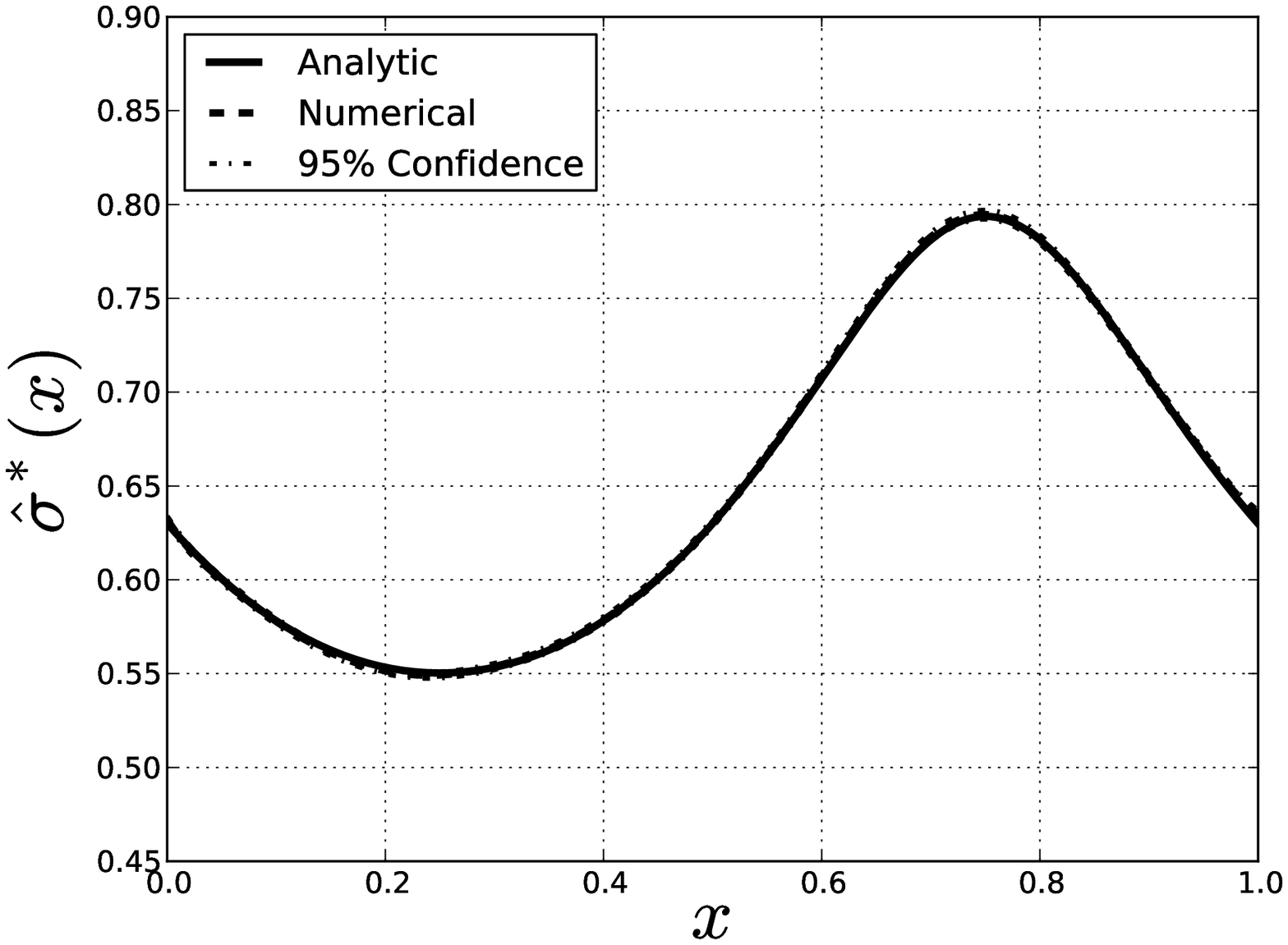}}
\caption{Optimal solutions obtained with increasing number of Monte Carlo samples. The gradient information used to obtain $\hat{\sigma}^*$ is computed using the sensitivity of the K-L expansion.}
\label{fig:sigma_opt}
\end{figure}

\begin{figure}[htbp]
\centering
\subfloat[$N = 10^2$]{\includegraphics[width=0.49\textwidth]{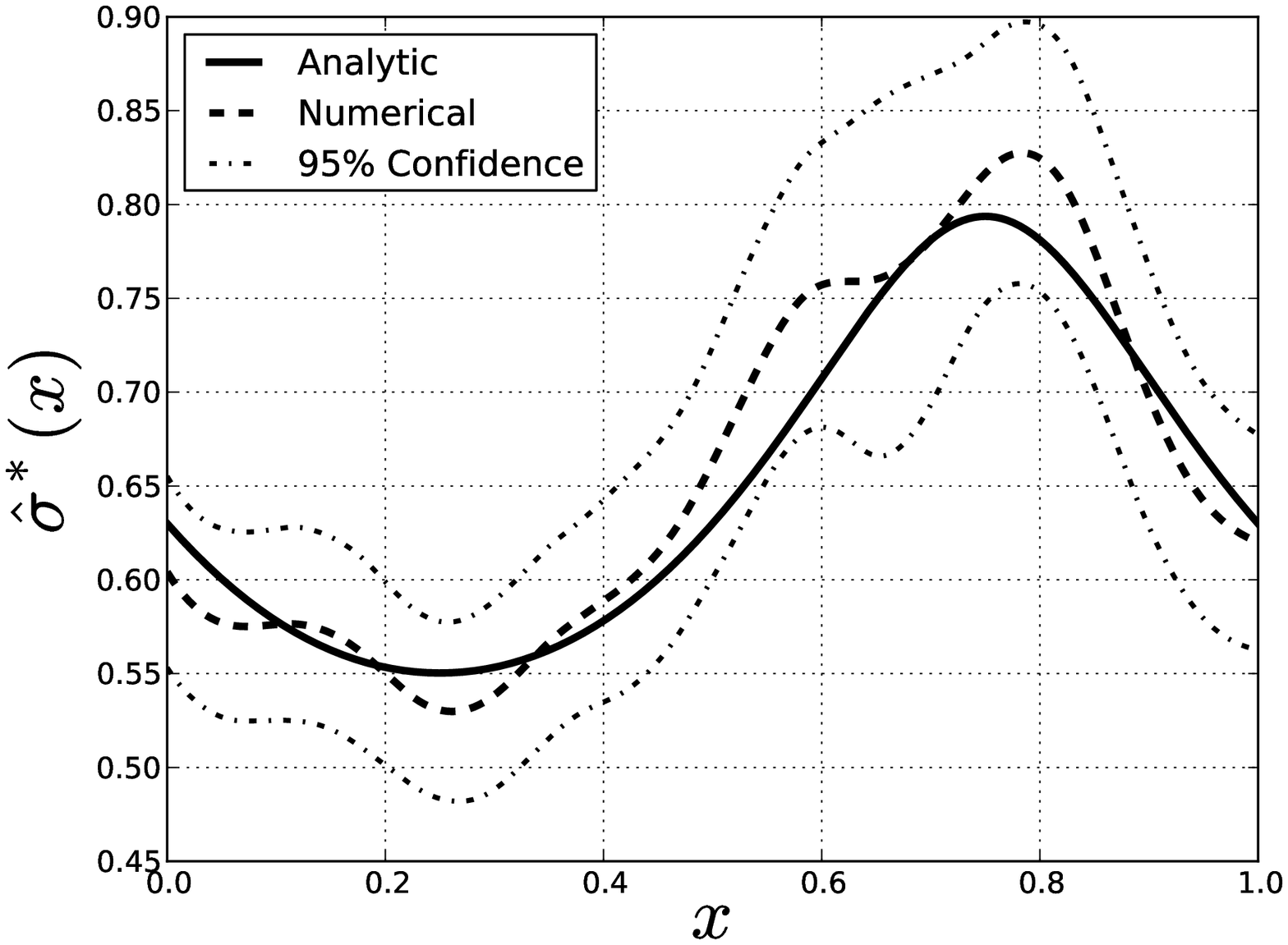}}
\hspace{0.1cm}
\subfloat[$N = 10^3$]{\includegraphics[width=0.49\textwidth]{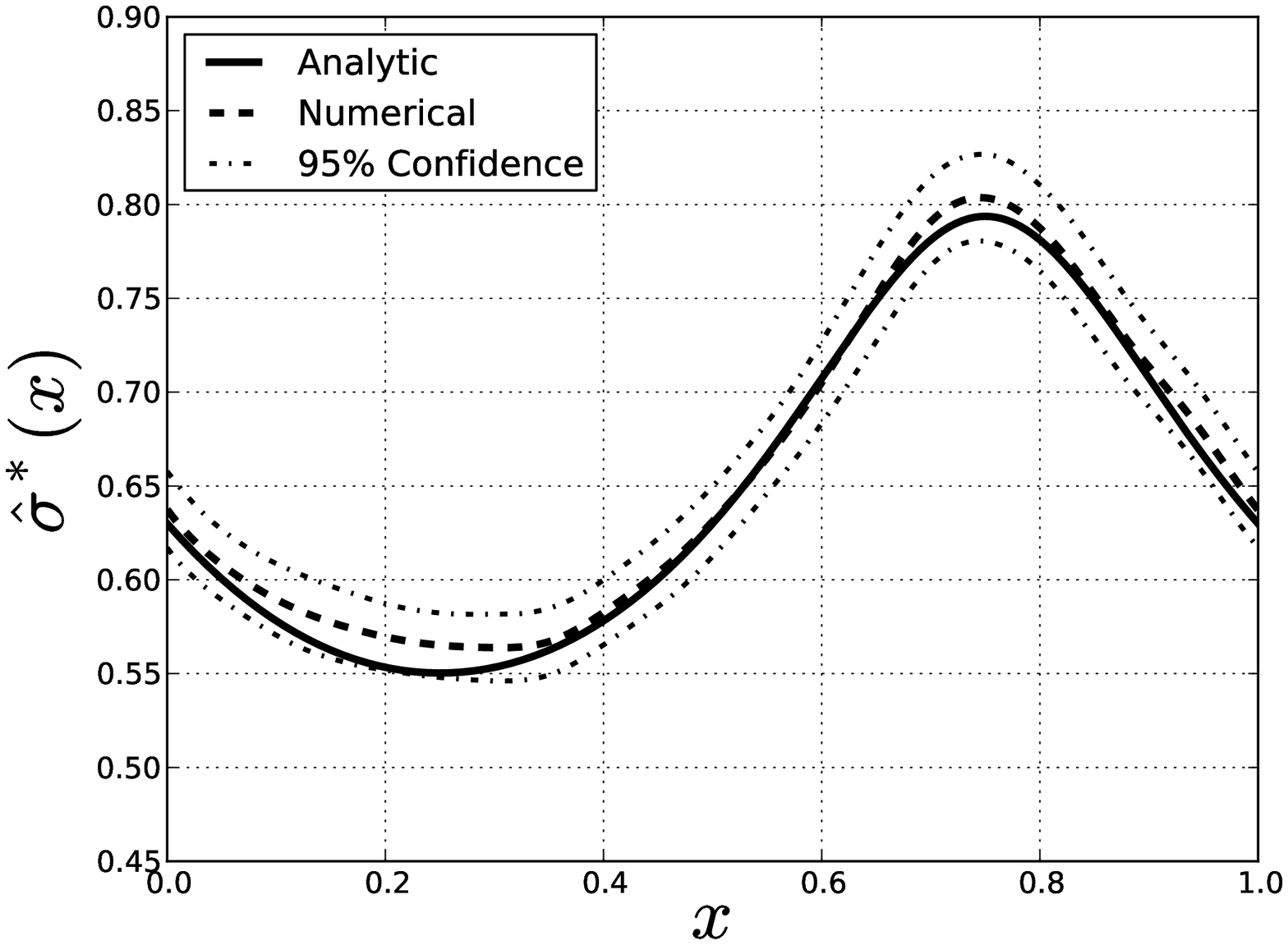}}
\newline
\subfloat[$N = 10^4$]{\includegraphics[width=0.49\textwidth]{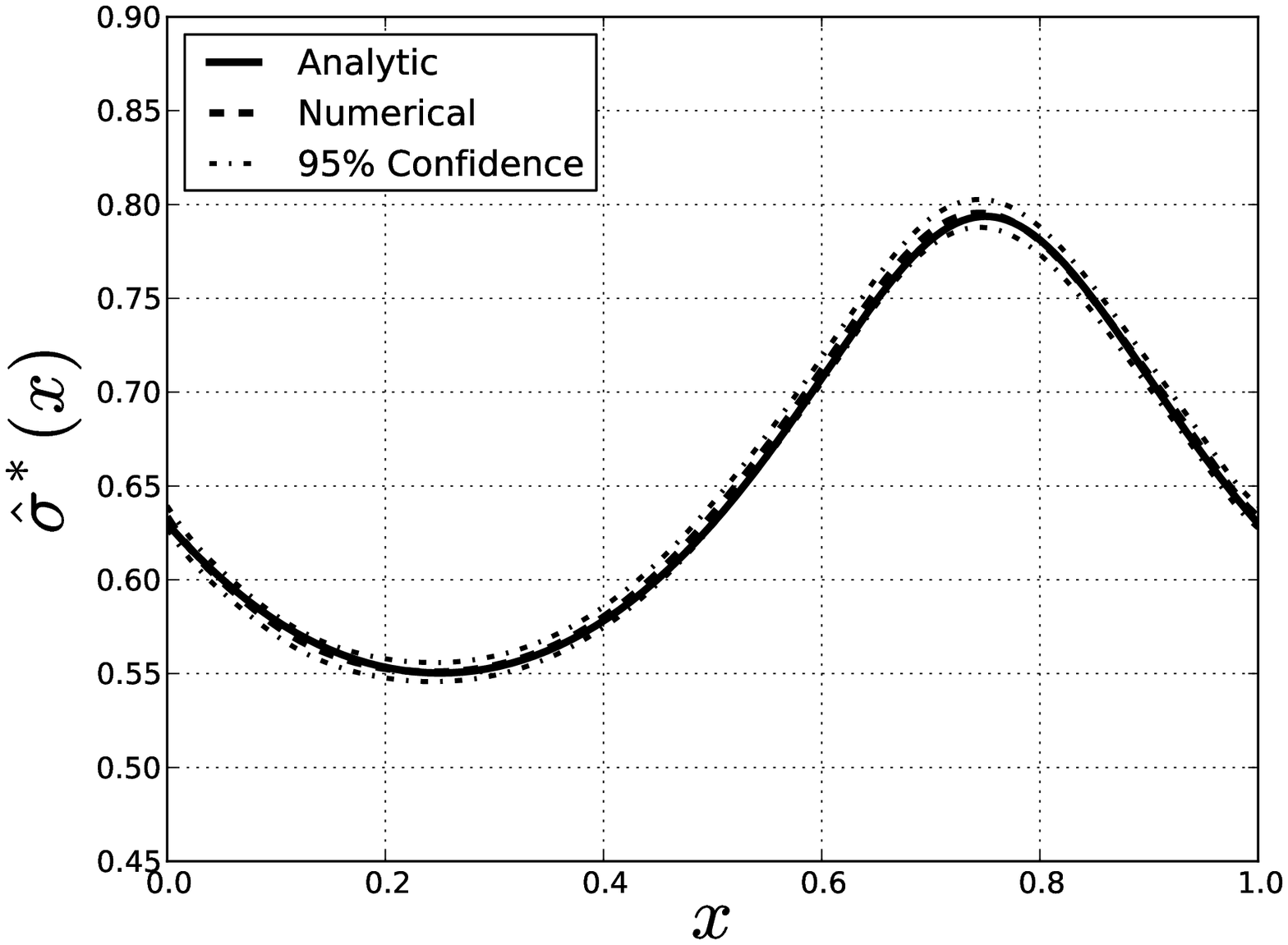}}
\hspace{0.1cm}
\subfloat[$N = 10^5$]{\includegraphics[width=0.49\textwidth]{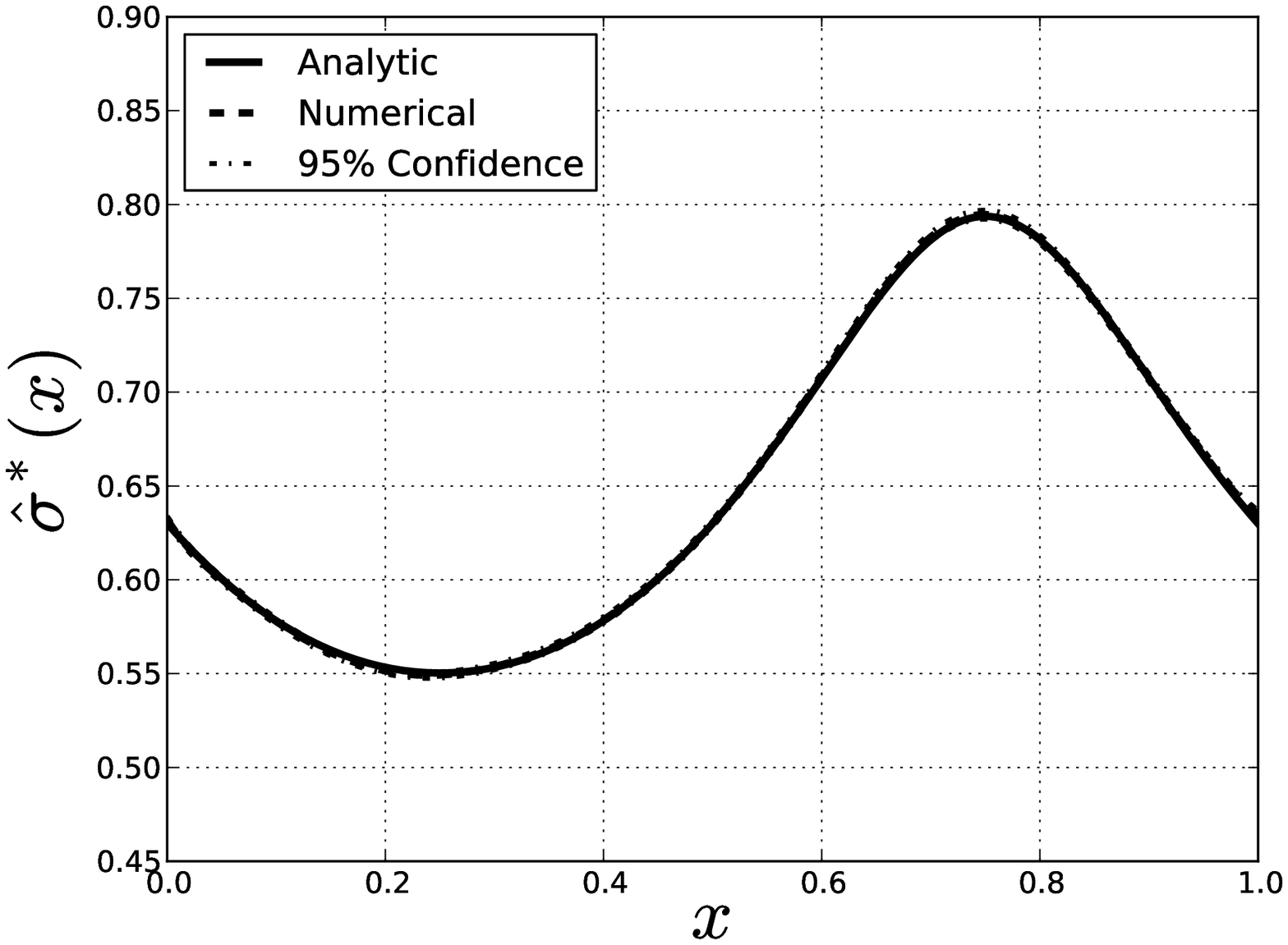}}
\caption{Optimal solutions obtained with increasing number of Monte Carlo samples. The gradient information used to obtain $\hat{\sigma}^*$ is computed using a scaled unit-variance random field.}
\label{fig:sigma_opt_fast}
\end{figure}

Figures \ref{fig:sigma_opt} and \ref{fig:sigma_opt_fast} show optimal solutions obtained using each approach. The shaded blue 95\% confidence region is computed by estimating the Hessian matrix $\mathbf{B}$ and covariance $\boldsymbol \Sigma$ using the Monte Carlo samples used to compute the optimal solution:
\begin{equation}
\hat{\mathbf{B}} = \frac{1}{N}\sum_{n=1}^{N} \nabla^2 F(\hat{\sigma}^*),
\end{equation}
\begin{equation}
\hat{\boldsymbol \Sigma} = \frac{1}{N}\sum_{n=1}^{N} \nabla F(\hat{\sigma}^*) \nabla F(\hat{\sigma}^*)^\intercal.
\end{equation}
The standard error of the optimal solution is then $\varepsilon_N = [\text{diag}(\hat{\mathbf{B}}^{-1}\hat{\boldsymbol\Sigma}\hat{\mathbf{B}}^{-1})/N]^{1/2}$. The plots show that the true optimal solution is largely within the 95\% confidence region for each approximate solution. Qualitatively, for a given number of Monte Carlo samples, the solutions obtained using either sensitivity approach are very similar. 

\begin{figure}[htbp]
\centering
\subfloat[$N = 10^2$]{\includegraphics[width=0.49\textwidth]{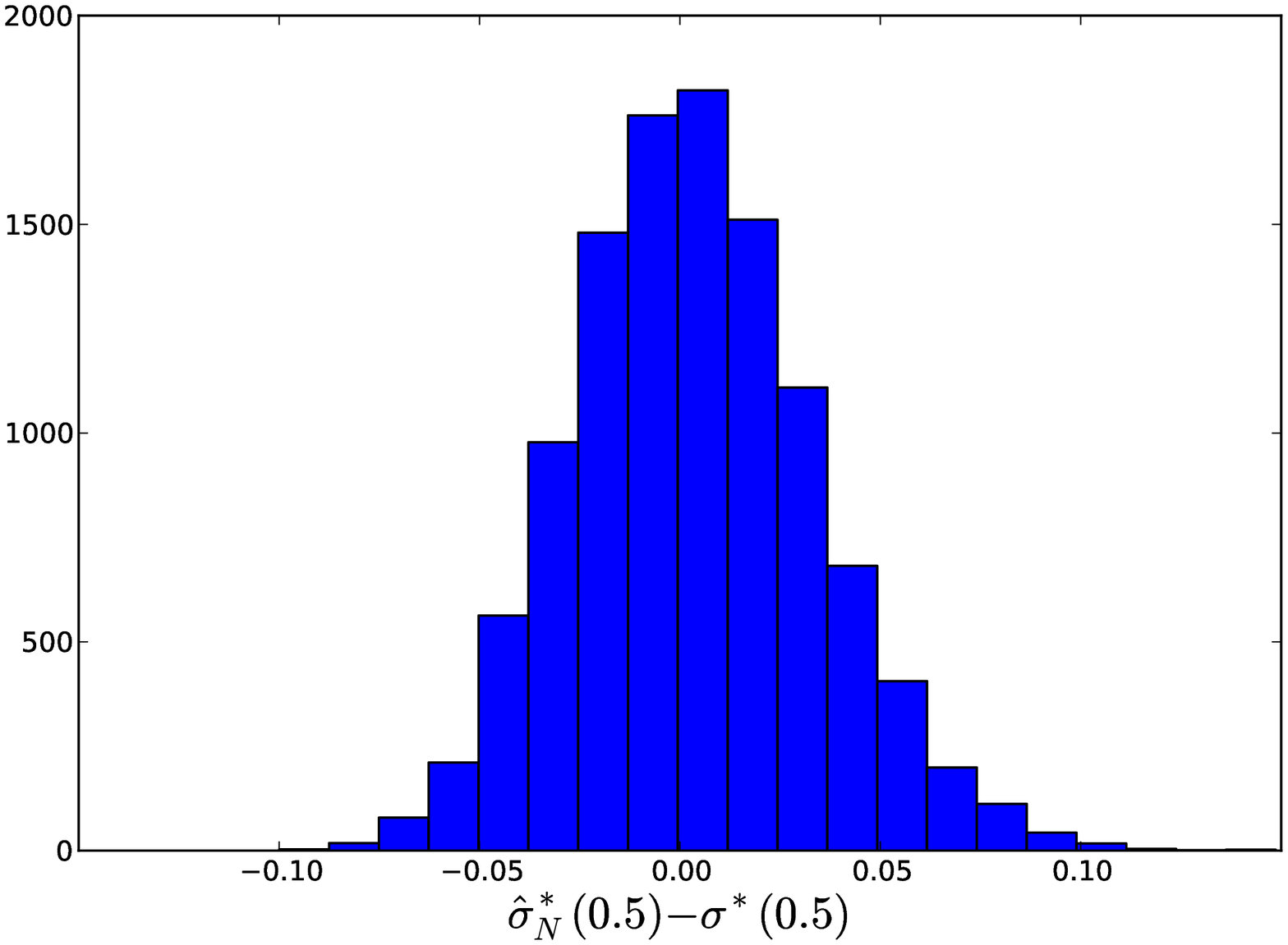}}
\hspace{0.1cm}
\subfloat[$N = 10^3$]{\includegraphics[width=0.49\textwidth]{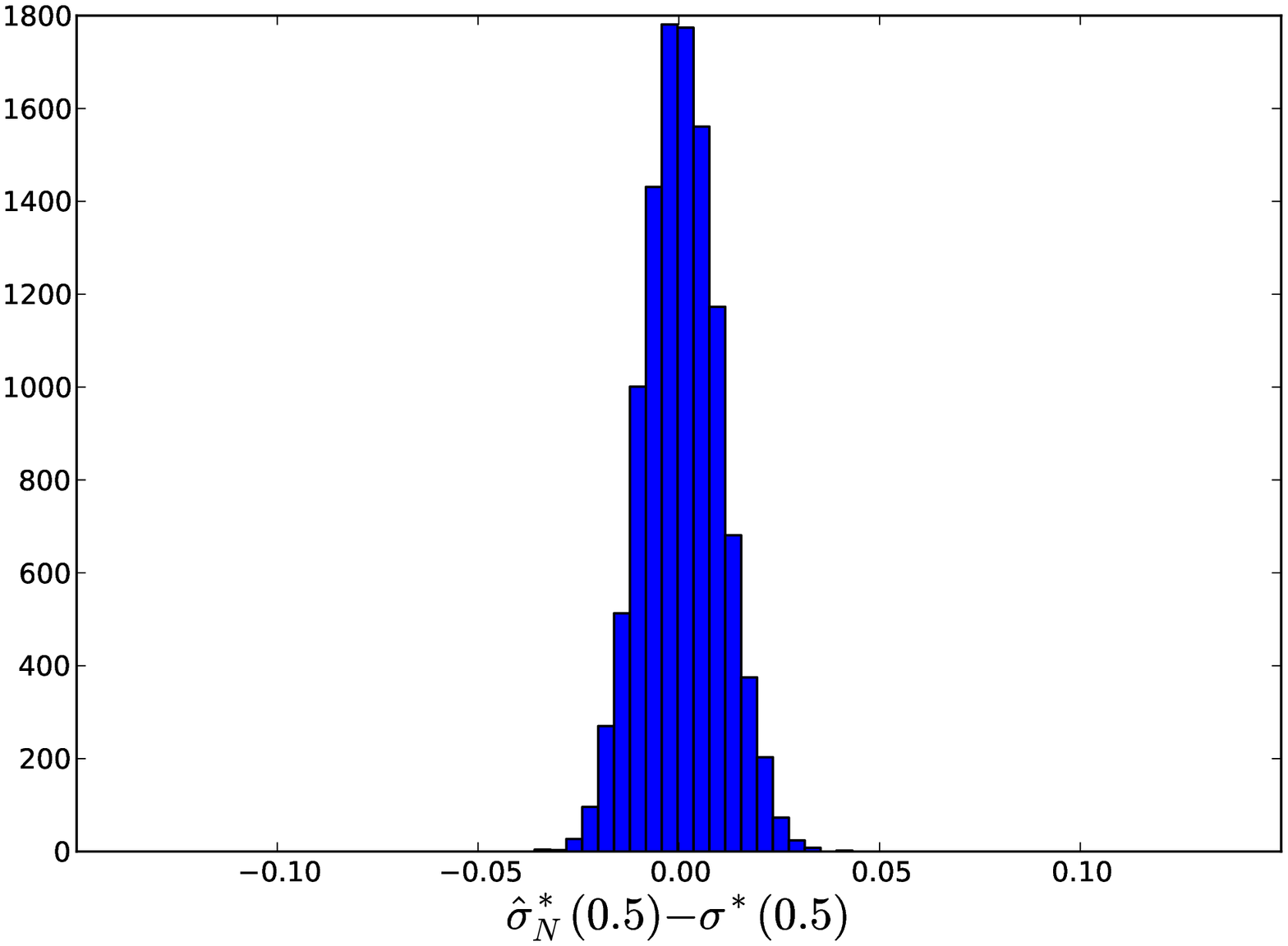}}
\newline
\subfloat[$N = 10^4$]{\includegraphics[width=0.49\textwidth]{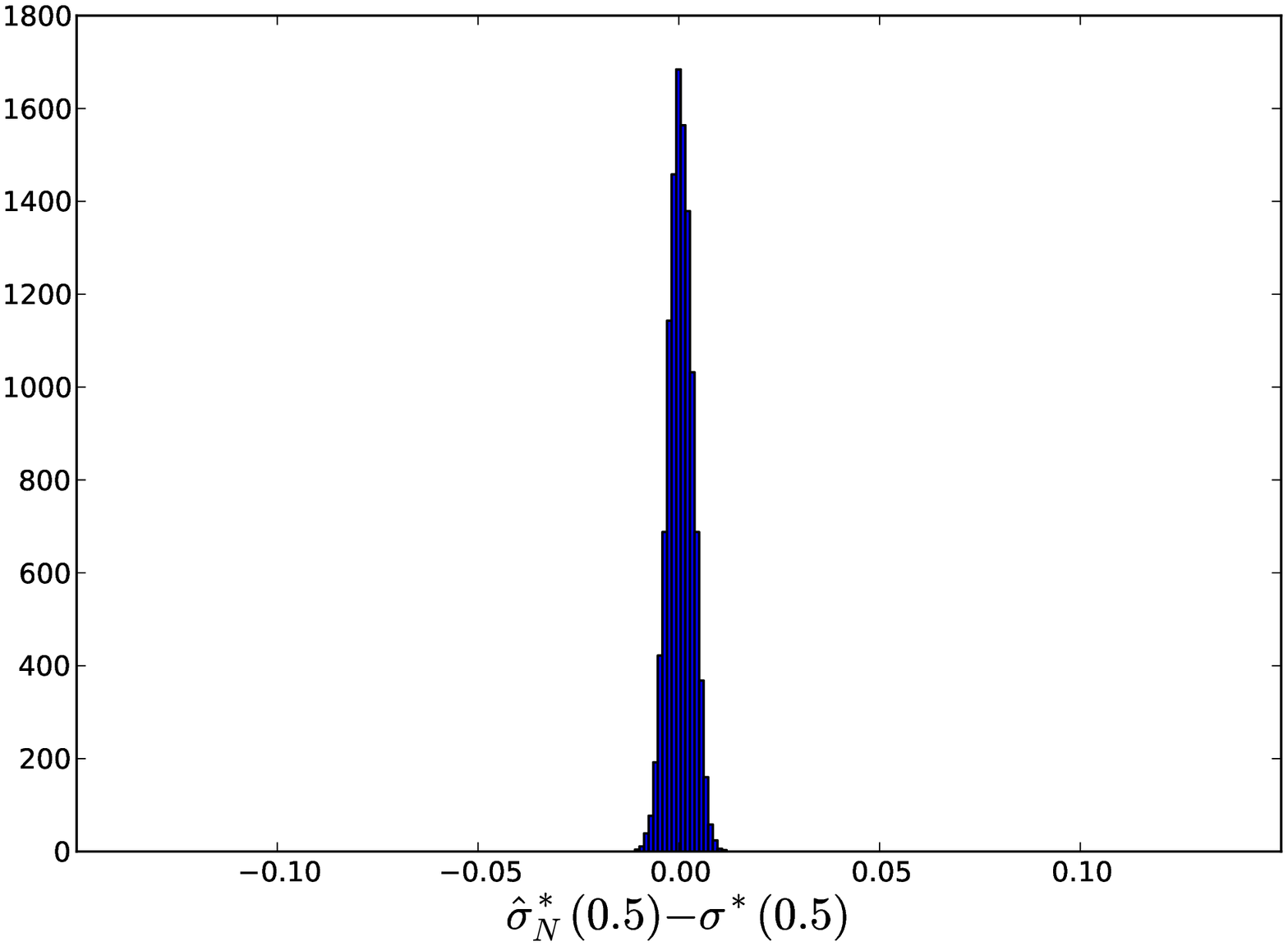}}
\hspace{0.1cm}
\subfloat[$N = 10^5$]{\includegraphics[width=0.49\textwidth]{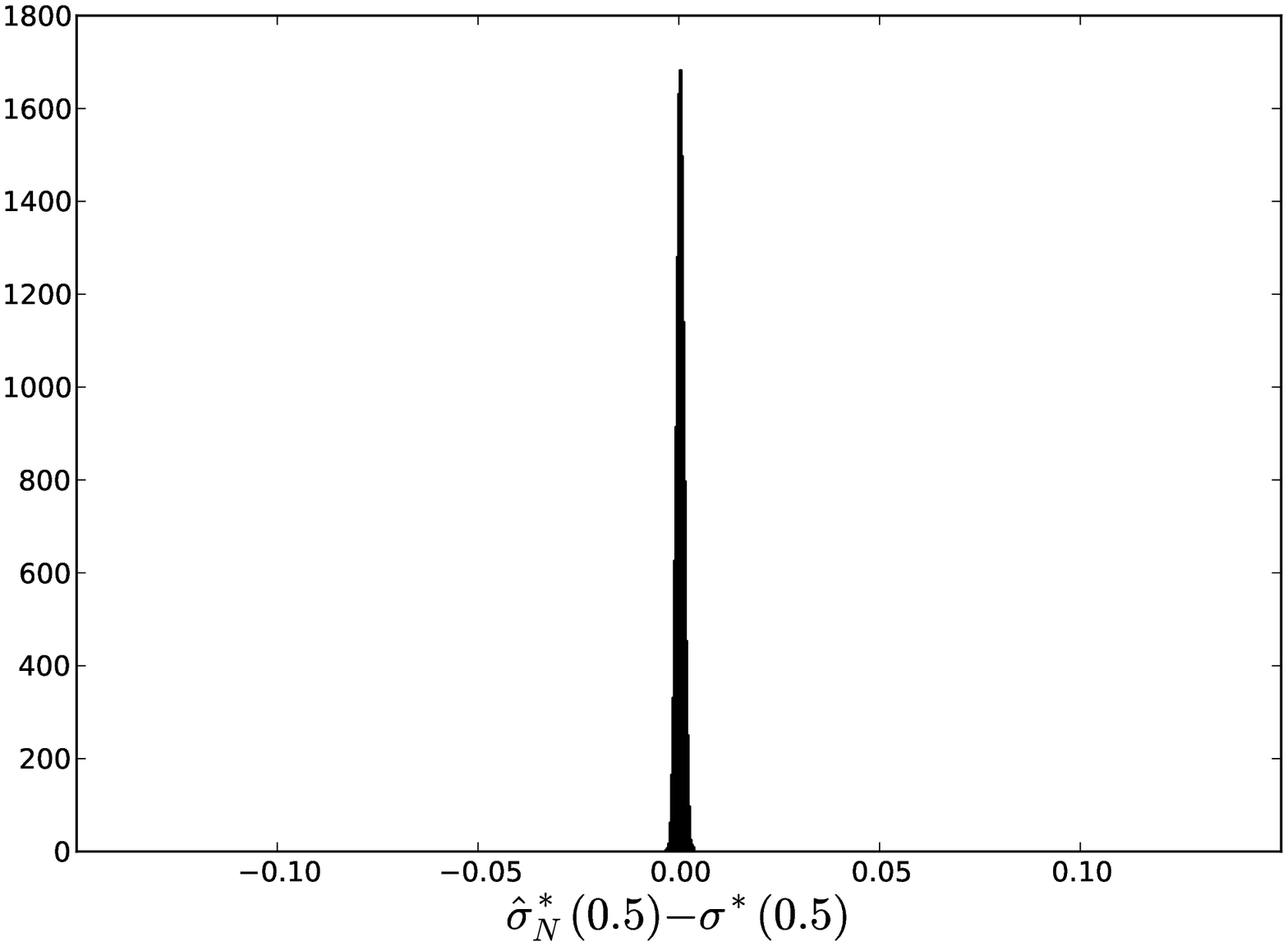}}
\caption{Histograms of the error at the center of the domain $\hat{\sigma}_N^*(0.5) - \sigma^*(0.5)$ for increasing number of Monte Carlo samples.}
\label{fig:hist_mid}
\end{figure}

\begin{figure}[htbp]
\subfloat[]{\includegraphics[width=0.5\textwidth]{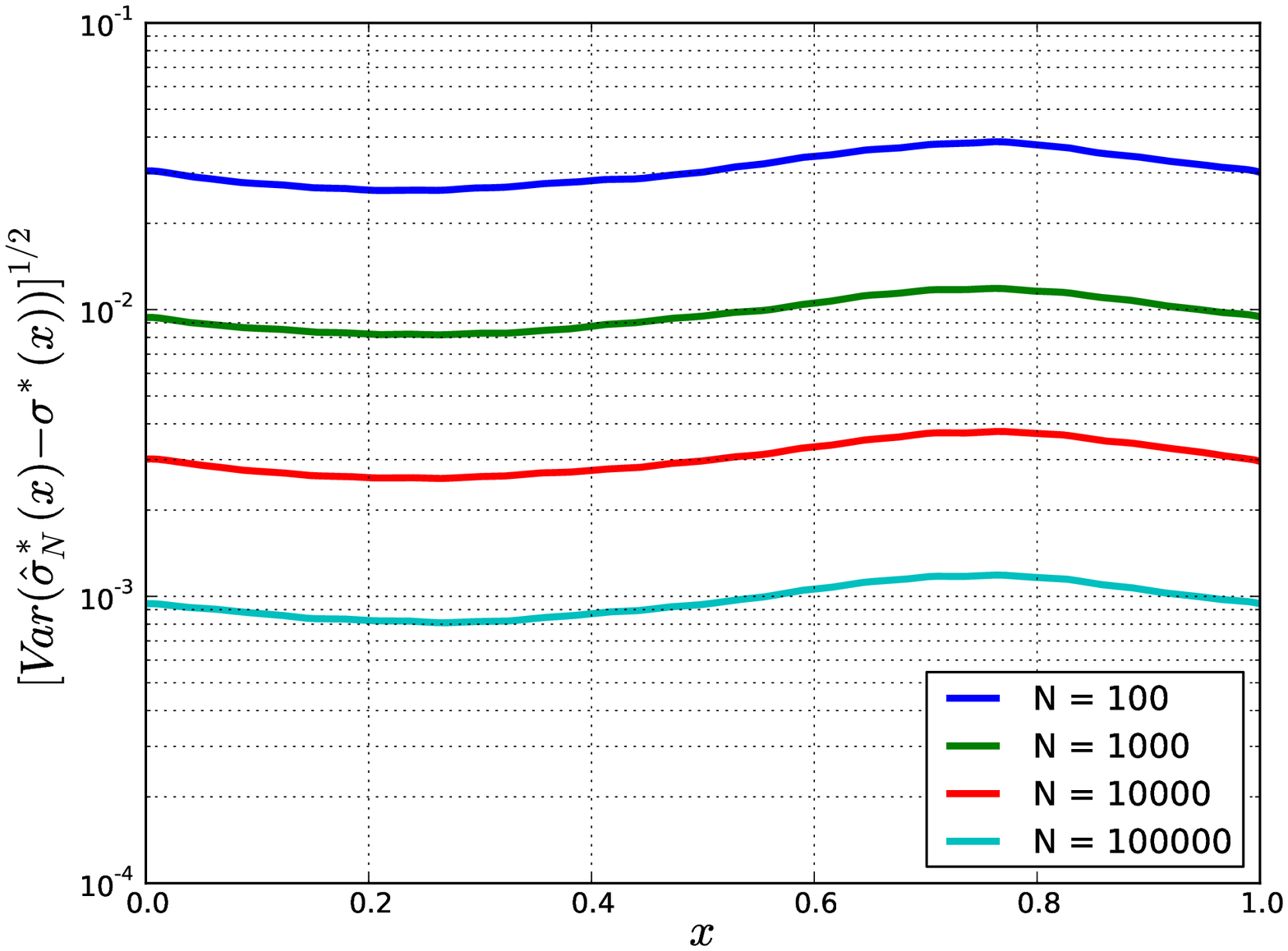}}
\hspace{0.1cm}
\subfloat[]{\includegraphics[width=0.5\textwidth]{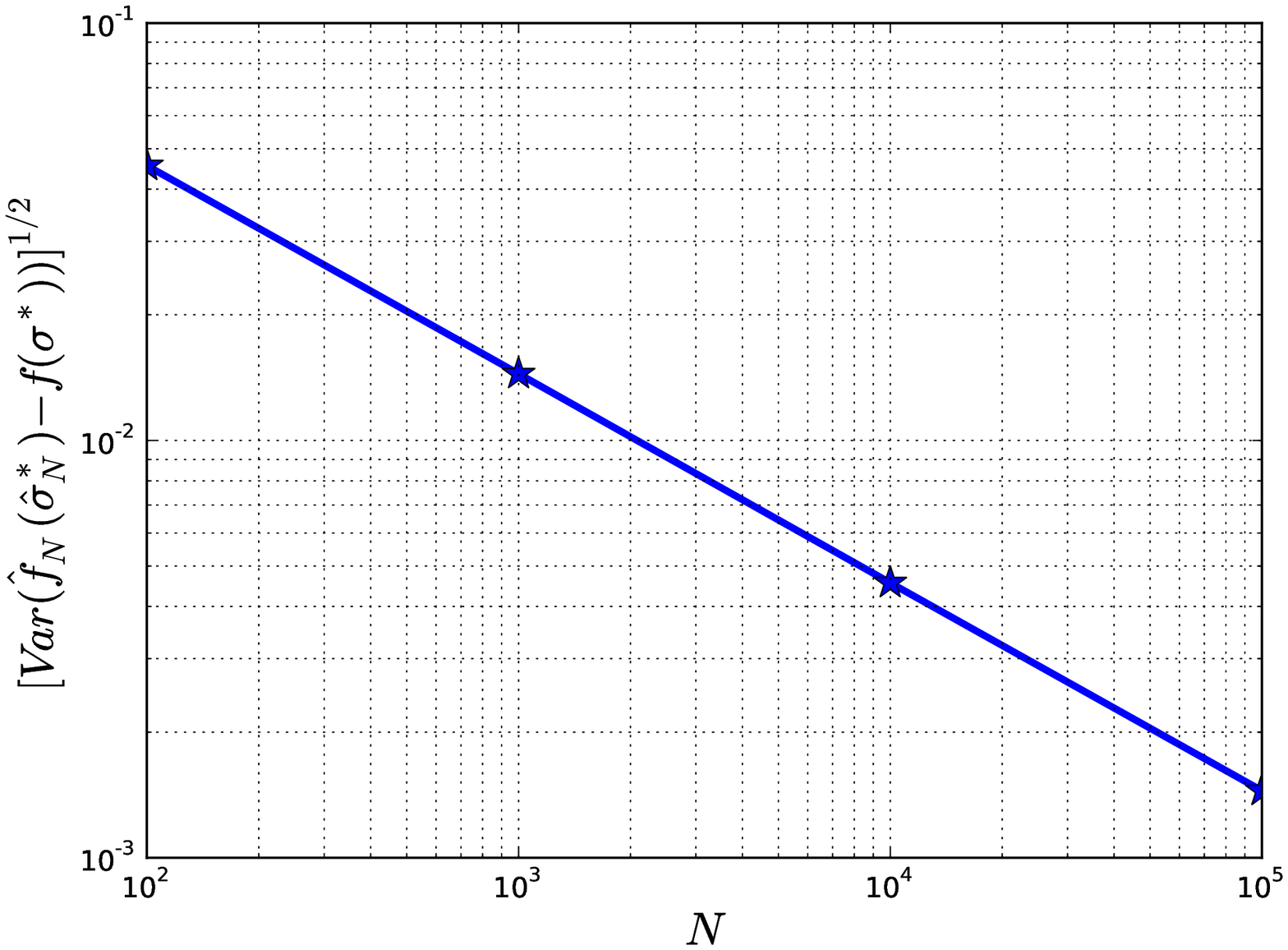}}
\caption{Standard deviations of the optimal solution (a) and optimal solution value (b).}
\label{fig:opt_conv}
\end{figure}

To further illustrate the convergence of the SAA optimal solution to the true optimal solution, we conduct $M = 10^4$ independent optimization runs for different values of $N$. This allows us to examine the distribution of the approximate optimal solution. Since the computational cost of using a scaled unit-variance random field is lower, we use this method to perform each optimization. Figure \ref{fig:hist_mid} shows histograms of the error of the SAA optimal solution evaluated at the center of the domain, i.e. $\hat{\sigma}_N^*(0.5) - \sigma^*(0.5)$, for various values of $N$. As expected, the histograms closely resemble Gaussian distributions with standard deviation proportional to $N^{-1/2}$. Figure \ref{fig:opt_conv} illustrates the convergence of the entire optimal solution and optimal value as $N$ is increased. The standard deviation of the optimal solution error $\hat{\sigma}_N^*(x) - \sigma^*(x)$ is plotted on the left, and the standard deviation of the optimal value error $f(\hat{\bp}_N^*) - f(\bp^*)$ is plotted on the right. We note that both converge like $N^{1/2}$: increasing the number of Monte Carlo samples by a factor of 100 gains a one decimal improvement in solution accuracy.

\section{Application: compressor blade tolerance optimization}

We now consider an application with engineering relevance: manufacturing tolerance optimization. Specifically, we consider a two-dimensional gas turbine compressor blade that is subject to geometric variability, and determine tolerances for this variability that provides the greatest performance benefit.

\subsection{Manufacturing error and tolerance models}

Previous studies of geometric variability in compressor blades has indicated that the discrepancy between manufactured blade geometries and the design intent geometry can be accurately modeled as a Gaussian random field\cite{garzon_thesis, sinha_2008}. In this context, the random field $e(x,\samp)$ represents the error between the manufactured surface and the nominal surface in the normal direction at the point $x$ on the nominal blade surface. The mean of the manufacturing error is assumed to be zero everywhere, i.e. $\bar{e}(x) = 0$. 

Manufacturing deviations tend to negatively impact the mean performance of compressor blades. We quantify the performance in terms of the total pressure loss coefficient, denoted by $\loss$, which measure the thermodynamic losses generated by a compressor. The mean total pressure loss coefficient tends to increase as the level of variability, i.e. the variance of the random field $e(x,\samp)$, is increased\cite{garzon_thesis}. It is possible to reduce this detrimental impact by specifying stricter manufacturing tolerances, thereby reducing the variance of the surface variations. 
% The relation between tolerances and variance is illustrated in figure~\ref{fig:tol_var}. The leading edge of the nominal blade geometry is shown in black. The $2\sigma$ bounds on the geometric variability are shown in dashed red lines, and sample realizations of manufactured blade geometries are plotted in color. The left figure demonstrates loose tolerances, where larger variations are allowed as compared to the case where strict tolerances are imposed, shown on the right. In each case, the variance $\sigma^2(x)$ of the random field is also non-stationary, decreasing near the leading edge in both cases.

%\begin{figure}[htbp]
%\centering
%\vspace{0.2cm}
%\subfloat[]{\includegraphics[width=0.45\textwidth]{loose_tol.png}}
%\hspace{1cm}
%\subfloat[]{\includegraphics[width=0.45\textwidth]{strict_tol.png}}
%\caption{Illustration of loose (a) and strict (b) tolerances.}
%\label{fig:tol_var}
%\end{figure}

To represent the standard deviation field $\sigma(x)$ over the surface of the blade, we use the same cubic B-spline basis introduced previously. The knot placement is chosen to enrich the basis near the leading edge, since previous studies of the impact of geometric variability on compressor performance have shown that most of the increase in loss results from imperfections near the leading edge\cite{garzon_2003}. A total of $N_\sigma = 31$ basis functions were used to parameterize the standard deviation.

\subsection{Flow solver}

All flow solutions are computed using the MISES (Multiple blade Interacting Streamtube Euler Solver) \cite{drela_giles_1987} turbomachinery analysis code. The boundary layer and wake regions are modeled using an integral boundary layer equation formulation describing the evolution of the integral momentum and kinetic energy shape parameter. In the inviscid regions of the flow field, the steady state Euler equations are discretized over a streamline conforming grid. Transition models are included to predict the onset of turbulent flow in the boundary layer.

A convenient feature of MISES is its solution speed. A typical flow solution requires 10-20 Newton Rhapson iterations to converge, which can be performed in a few seconds. Moreover, MISES offers the option to reconverge a flow solution after perturbing the airfoil geometry. Since the perturbations in the geometry introduced by manufacturing variability are small, the flow field corresponding to blades with manufacturing variability can be reconverged very quickly from the flow field computed for the nominal geometry. This offsets some of the computational cost associated with using the standard Monte Carlo method to propagate uncertainty.

%The nominal blade surface is parameterized using a combination of Chebyshev polynomial modes, acting in the direction normal to the blade surface, and a stagger angle mode to rotate the blade around the leading edge. Five Chebyshev modes were used, giving a total of $N_d = 6$ nominal geometry design parameters. The vector of design vectors is represented by $\bd \in \mathbb{R}^{N_d}$. To ensure the aeromechanical response of the blade is not significantly altered, the blade thickness is constrained by applying equal and opposite perturbations to corresponding points on the pressure and suction surfaces. This choice of geometric design modes provides an orthonormal basis onto which to expand modifications to the camber of the airfoil, as well as providing increased design resolution near the leading and trailing edges.

\subsection{Optimization statement}

We seek to optimize the manufacturing tolerances to reduce the detrimental impact of manufacturing variations. To do this, we first define the variability metric $V$, which measures the total level of manufacturing variations over the entire blade surface:
\begin{equation}
V(\bs) = \int_X \sigma(x)\,dx.
\end{equation}
Here $\bs \in \mathbb{R}^{N_\sigma}$ parameterizes the standard deviation $\sigma(x)$. Specifying stricter tolerances (decreasing $V$) incurs higher manufacturing costs. To constrain this cost, we constrain the variability metric to a specific value $V_b$, representing the strictest tolerances deemed acceptable by the manufacturer. The standard deviation of the manufacturing variability is constrained from above to ensure the optimizer does not trade increases in variability in regions of low sensitivity for excessive decreases in variability in regions of high sensitivity. The resulting optimization problem for the optimal tolerances is given below.
\begin{equation}
\begin{aligned}
& \bs^* =
& & \underset{\bs}{\arg \min}
& & \E[\loss(\bs)] \\
& & & \ \ \ \ \text{s.t.}
& & V(\bs) = V_b \\
& & & & & \sigma(x) \leq \sigma_\text{max}
\end{aligned}
\label{eq:comp_opt}
\end{equation}
To solve (\ref{eq:comp_opt}) numerically, the SAA method is used and all objective and constraint functions are replaced by their Monte Carlo estimates. The resulting nonlinear optimization problem is solved using SQP. The gradient of the objective and constraints is computed using the pathwise approach described previously. The shape sensitivities are evaluated using second-order accurate finite differences.

\subsection{Numerical results}

We apply the proposed method to optimize the tolerances of a two-dimensional fan exit stator cascade. In the absence of geometric variability, the loss coefficient is $\loss = 2.22 \times 10^{-2}$. Manufacturing variations are prescribed in the form of a Gaussian random field with standard deviation $8.0 \times 10^{-4}$ (non-dimensionalized by the blade chord). The covariance function of the random field is the same squared exponential function described earlier, with a correlation length $L$ that is reduced near the leading edge of the blade to reflect the manufacturing variations observed in measured blades. The mean loss coefficient of the blade in the presence of manufacturing variability is $\E[\loss] = 2.29 \times 10^{-2}$, roughly 3\% higher than the loss of the design intent geometry. 

A total of 75 SQP iterations were required to obtain the optimal solution. Each SQP iteration requires evaluating one evaluation of the gradient of the objective and constraint functions, as well as a number of evaluations of the objective and constraint functions to perform a linesearch. This resulted in 120 Monte Carlo simulations, each comprised of $N = 500$ flow solutions. Thus, the parallelizability of the Monte Carlo method and the speed of the MISES code had great benefit.

\begin{figure}[htbp]
\centering
\vspace{0.5cm}
\includegraphics[width=0.6\textwidth]{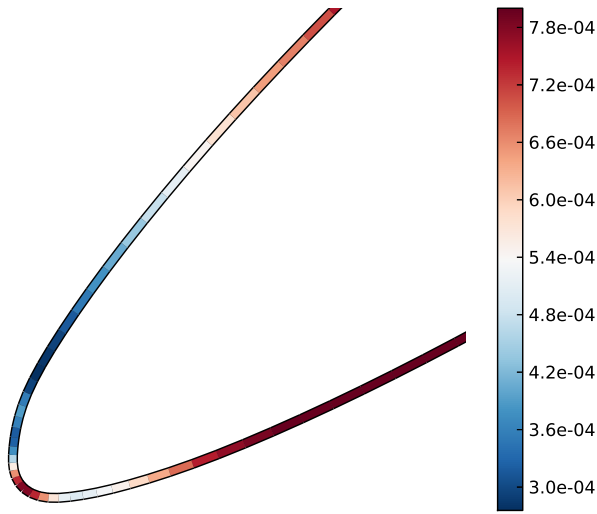}
\caption{Optimal distribution of the standard deviation $\sigma(x)$.}
\label{fig:opt_results}
\end{figure}

The total allowable variability $V_b$ was constrained to be 98\% of the baseline level of variability. The optimized tolerances are shown in Figure \ref{fig:opt_results}. We only show the standard deviation near the leading edge of the blade, since the optimal value over the rest of the blade was equal to the baseline value of $8.0 \times 10^{-4}$. We observe that the greatest reduction in variability is specified on the upper surface of the blade. The optimized tolerances reduce the mean loss coefficient to $\E[\loss] = 2.23 \times 10^{-2}$, roughly 0.5\% higher than the loss of the design intent geometry. For a very small decrease in the level of manufacturing variability, a significant increase in the mean performance is realized, demonstrating the efficacy of the proposed approach.

\section{Summary and conclusions}

Considerable research has been conducted in the area of design under uncertainty, bringing together the fields of uncertainty quantification and optimization. Optimization of the uncertainty itself has received considerably less attention. This paper has presented an approach for optimizing the mean and covariance of Gaussian random fields to achieve a desired statistical performance. The novel sensitivity analysis presented here allows for gradient-based algorithms to be leveraged when performing these optimizations. 

The approach presented in this paper can be applied when the mean and covariance functions depend explicitly on some set of parameters. We have presented the example of tolerance optimization, where the level of variability is a design variable. Another example arises from optimizing measurement locations in a Gaussian random field, where, conditioned on the measurements, the covariance depends explicitly on the measurement location.
Future improvements to the proposed optimization framework would incorporate adjoint sensitivity information when considering PDE-constrained problems. This would reduce the computational cost of estimating gradients when the number of design parameters is large with respect to the number of objectives and constraints, which is common in engineering optimization.

%\Appendix
%\section{The use of appendices}
%The \verb|\appendix| command may be used before the final sections
%of a paper to designate them as appendices. Once \verb|\appendix|
%is called, all subsequent sections will appear as 
%
%\appendix
%\section{Title of appendix} Each one will be sequentially lettered
%instead of numbered. Theorem-like environments, subsections,
%and equations will also have the section number changed to a letter. 

\clearpage
\bibliographystyle{siam}
\bibliography{main_sisc}

\end{document}